\documentclass{article}%
\usepackage{amsmath}
\usepackage{amsfonts}
\usepackage{amssymb}
\usepackage{graphicx}%
\providecommand{\U}[1]{\protect\rule{.1in}{.1in}}

\begin{document}

\title{Combining Evidence}
\author{Michael Evans and Yang Guan Jian Guo\\Department of Statistical Sciences, University of Toronto}
\date{}
\maketitle

\begin{center}
\textbf{Abstract}
\end{center}

The problem of combining the evidence concerning an unknown, contained in each
of $k$ Bayesian inference bases, is discussed. This can be considered as a
generalization of the problem of pooling $k$ priors to determine a consensus
prior. The linear opinion pool of Stone (1961) is seen to have the most
appropriate properties for this role. In particular, linear pooling preserves
a consensus with respect to the evidence and other rules do not. While linear
pooling does not preserve prior independence, it is shown that it still
behaves appropriately with respect to the expression of evidence in such a
context. For the general problem of combining evidence, Jeffrey
conditionalization plays a key role.\smallskip

\noindent\textbf{Keywords and phrases}: \textit{combining priors, statistical
evidence, preserving consensus, Jeffrey conditionalization, ancillarity.}

\section{Introduction}

Suppose that $k$ different experts choose models and priors for a statistical
analysis concerning a common quantity of interest $\Psi$ which is a parameter
or a future value. A problem then arises as to how the resulting statistical
analyses should be combined so that the inferences presented can serve as a
consensus inference. If all the models are the same, then this is the
well-known problem of combining priors and this is covered by our discussion
here. Even for the problem of combining priors, however, a somewhat different
point-of-view is taken. A particular measure of evidence is adopted, as
discussed in Section 3, such that the data set, sampling model and prior leads
to either evidence in favor of or against each possible value of $\Psi$. The
purpose then is to determine a consensus on what the evidence indicates. So,
even for the combining priors problem, the motivation here is combining the
evidence rather than combining prior beliefs. Since the primary goal of a
statistical analysis is to express what the evidence says, this seems
appropriate. Also, it is perfectly reasonable that some analyses express
evidence against while others express evidence in favor but the combined
expression of the evidence is one way or the other, see Section 2.

Before discussing the combination approach, however, it is necessary to be
more precise about the problem and distinguish between somewhat different
contexts where the problem can arise. It will be supposed here that $\Psi$ is
a parameter of interest but prediction problems are easily handled by a slight
modification, see Example 3. Let $\mathcal{M}=\{f_{\theta}:\theta\in\Theta\}$
denote a generic statistical model and $\psi=\Psi(\theta),$ where $\Psi
:\Theta\rightarrow\Psi$ is onto and to save notation the function and its
range have the same symbol.\smallskip

\noindent\textbf{Context I. }Suppose there is a single statistical model
$\mathcal{M}$ for the data $x$ and $k$ distinct priors $\pi_{i}$ so there are
$k$ inference bases $\mathcal{I}_{i}=(x,\mathcal{M},\pi_{i})$ for
$i=1,\ldots,k.$ It is assumed that the conditional priors $\pi_{i}%
(\cdot\,|\,\psi)$ on the nuisance parameters are all the same, as is satisfied
when $\Psi(\theta)=\theta.$ This situation arises when there is a group of
analysts who agree on $\mathcal{M}$ and perhaps use a default prior for the
nuisance parameters, while each member puts forward a prior for $\Psi
.$\smallskip

\noindent\textbf{Context II. }Suppose there are $k$ data sets, models, and
priors as given by the inference bases $\mathcal{I}_{i}=(x_{i},\mathcal{M}%
_{i},\pi_{i})$ for $i\in\{1,\ldots,k\}$ and there is a common characteristic
of interest $\psi=\Psi(\theta_{i})$ with the true value of $\psi$ being the
same for each model, as will occur when $\psi$ corresponds to some real world
quantity.\smallskip

\noindent Note that $\psi$ references some real-world quantity so in Context
II the set of possible values and its true value is the same for each model
even though formally the function $\Psi$ may differ between models but we
suppress this in the notation.

It is a necessary part of any statistical analysis that a model be checked to
see whether or not it is contradicted by the data, namely, determining if it
is the case that the data lies in the tails of each distribution in the model.
So in any situation where there is a lack of model fit, it is necessary to
modify that component of the inference base. Similarly, each prior needs to be
checked for prior-data conflict, namely, is there an indication that the true
value lies in the tails of the prior, see Evans and Moshonov (2006), Nott et
al. (2020). If such a conflict is found, then the prior needs to be modified,
see Evans and Jang (2011). It is assumed hereafter that all the models and
priors have passed such checks. A salutary effect of a lack of prior-data
conflict is that it rules out the possibility of trying to combine priors
which have little overlap in terms of where they place their mass.

Given an inference base $\mathcal{I}=(x,\mathcal{M},\pi)$ and interest in
$\psi=\Psi(\theta),$ a Bayesian analysis has a consistency property. In
particular, this inference base is equivalent to the inference base
$\mathcal{I}=(x,\mathcal{M}_{\Psi},\pi_{\Psi})$ for inference about $\psi,$
where $\pi_{\Psi}$ is the marginal prior on $\psi$ and $\mathcal{M}_{\Psi
}=\{m_{\psi}:\psi\in\Psi\}$ with $m_{\psi}(x)=E_{\Pi(\cdot\,|\,\psi
)}(f_{\theta}(x)),$ the prior predictive density of the data obtained by
integrating out the nuisance parameters via the conditional prior $\Pi
(\cdot\,|\,\psi)$ for $\theta$ given $\psi=\Psi(\theta).$ So, for example, the
posterior $\pi_{\Psi}(\cdot\,|\,x)$ for $\psi$ obtained via these two
inference bases is the same and moreover the evidence about $\psi$ is also the
same. This result has implications for the combination strategy as it is
really the inference bases $\mathcal{I}_{i}=(x,\mathcal{M}_{\Psi},\pi_{i\Psi
})$ that are relevant in Context I and it is the inference bases
$\mathcal{I}_{i}=(x_{i},\mathcal{M}_{i\Psi},\pi_{i\Psi})$ that are relevant in
Context II, namely, nuisance parameters are always integrated out before combining.

Perhaps the simplest step away from Context I is when the data sets differ but
all the models are based on the same basic set of candidates for the true
probability measure and with the same conditional prior on the nuisance
parameters. In such a context it seems obviously correct to simply combine the
data sets and use the common model for the combined data set which places the
problem within Context I. The result would not necessarily be the same if the
data sets were not combined, so it is necessary that the following rule be
applied first to the set of inference bases in general.\smallskip

\noindent\textit{Combining inference bases rule}: all data sets that are
assumed to arise from the same set of basic distributions are combined so that
separate data sets are associated with truly distinct models.\smallskip

\noindent This rule ensures that any combination reflects true differences
among the beliefs concerning where the truth lies as there is agreement on the
other ingredients. It is assumed hereafter that this is applied before the
inference bases $\mathcal{I}_{i}$ are determined. Note too that, even if the
basic model is the same for each $i,$ when the conditional priors on the
nuisance parameters differ, then this is Context II.

In\ Section 2 a general family of rules for combining priors with given
weights is presented. In Section 3 the problem of combining evidence for
Context I is analyzed, with given weights for the respective priors, and the
linear pooling combination rule is seen to have appropriate properties with
respect to evidence. In\ Section 4 the problem of determining appropriate
weights is considered. In Section 5 the problem for Context II is discussed
and a proposal is made for a rule that generalizes the rule for Context I. The
rule for Context I possesses a natural consistency property as the combined
evidence is the same whether considered as a mixture of the evidence arising
from each inference base or obtained directly from the combined prior and the
corresponding posterior. In particular, it is Bayesian in this generalized
sense which differs from being externally Bayesian as discussed in Genest
(1984), see Section 3. This is not the case for Context II, however, because
of ambiguities in the definition of the likelihood, but Jeffrey
conditionalization provides a meaningful interpretation, at least when all the
inference bases contain the same data.

The problem of combining priors has an extensive literature. Winkler (1968) is
a basic reference and reviews can be found in Genest and Zidek (1986), Clemen
and Winkler (1999), O'Hagan et al. (2006) and French (2011). Farr et al.
(2018) is a significant recent application. Broadly speaking there are
mathematical approaches and behavioral approaches. The mathematical approach
provides a formal rule, as in Section 2, while the behavioral approach
provides methodology for a group of proposers to work towards a consensus
through mutual interaction. For example, Burgman et al. (2011) consider the
elicitation procedure where quantities concerning the object of interest are
elicited by each member of a group and then the average elicited values are
used to choose the prior. Albert et al. (2012) adopt a supra-Bayesian approach
where the data generated during the elicitation process is conditioned on in a
formal Bayesian analysis to choose a prior in a family on which an initial
prior has been placed. Yager and Alajlan (2015) present an iterative
methodology for a group of proposers to work towards a consensus prior based
upon each proposer seeing how far their proposal deviated from a current
grouped proposal. While the behavioral approach has a number of attractive
features, there are also reservations as indicated by Kahneman in Goodman (2021).

The focus in this paper is on presenting a consensus assessment of the
evidence via a combination of the evidence that each analyst obtains. In
particular, the priors $\pi_{i}$ need not arise via the same elicitation
procedure and the proposers may not be aware of other proposals although the
approach does not rule this out. Also, utility functions, necessary for
decisions, are not part of the development as these may indeed lead to
conflicts with what the evidence indicates and they are not generally
checkable against the data as with models and priors. The assessment of
statistical evidence as the primary driver of statistical methodology is a
theme that many authors have pursued, for example, Birnbaum(1962), Royall
(1997),\ Evans (2015), Gelman and O'Rourke (2017) and Vieland and Chang (2019).

\section{Combining Priors with Given Prior Weights}

Let $\alpha=(\alpha_{1},\ldots,\alpha_{k})\in S_{k}$ the $(k-1)$-dimensional
simplex for some $k\geq2$ and, for now, suppose that $\alpha$ is given. While
general combination rules could be considered, attention is restricted here to
the power means of densities%
\[
\pi_{t,\alpha}=\left\{
\begin{tabular}
[c]{ll}%
$c_{t}(\alpha,\pi_{\cdot})\{\sum_{i=1}^{k}\alpha_{i}\pi_{i}^{t}\}^{1/t}$ &
$\text{when }t\neq0,\pm\infty$\\
$c_{0}(\alpha,\pi_{\cdot})\exp\{%
{\textstyle\sum_{i=1}^{k}}
\alpha_{i}\log\pi_{i}\}$ & when $t=0$\\
$c_{-\infty}(\alpha,\pi_{\cdot})\min\left\{  \pi_{1},\ldots,\pi_{k}\right\}  $
& when $t=-\infty$\\
$c_{\infty}(\alpha,\pi_{\cdot})\max\left\{  \pi_{1},\ldots,\pi_{k}\right\}  $
& when $t=\infty$%
\end{tabular}
\ \right.
\]
where $\pi_{\cdot}=(\pi_{1},\ldots,\pi_{k})$ and, for any $\alpha$ and
sequence of nonnegative functions $g_{\cdot}=(g_{1},\ldots,g_{k})$ defined on
$\Theta,$ then $c_{t}(\alpha,g_{\cdot})$ is the relevant normalizing constant.
Note that $\pi_{-\infty,\alpha}$ and $\pi_{\infty,\alpha}$ do not depend on
$\alpha.$

For each $\theta$ the mean $\{\sum_{i=1}^{k}\alpha_{i}\pi_{i}^{t}%
(\theta)\}^{1/t}$ is nondecreasing in $t,$ see Steele (2004), and two of the
means are equal everywhere iff all priors are the same$.$ Since $c_{1}%
(\alpha,\pi_{\cdot})=1,$ this implies that $c_{t}(\alpha,\pi_{\cdot})$ is
finite for all $\alpha$ whenever $t\leq1.$ If $t>1$ is to be considered, then
it is necessary to check on the integrability of the mean so that a proper
prior is obtained and this will be assumed to hold whenever the case $t>1$ is
referenced. When $\Theta$ is finite, this is not an issue.

The following result characterizes how the posterior behaves in terms of a
combination of the individual posteriors. Let $m_{i}$ denote the $i$-th prior
predictive density based on prior $\pi_{i}$ and $m_{t,\alpha}$ denote the
prior predictive density obtained using the $\pi_{t,\alpha}$ prior.\smallskip

\noindent\textbf{Proposition} \textbf{1.} For Context I, the posterior based
on $\pi_{t,\alpha}$ equals%
\[
\pi_{t,\alpha}(\theta\,|\,x)=\left\{
\begin{tabular}
[c]{ll}%
$c_{t}(\alpha,m_{\cdot}(x)\pi_{\cdot}(\cdot\,|\,x))\{\sum_{i=1}^{k}\alpha
_{i}m_{i}^{t}(x)\pi_{i}^{t}(\theta\,|\,x)\}^{1/t}$ & $\text{when }t\neq0$\\
$c_{0}(\alpha,\pi_{\cdot}(\cdot\,|\,x))\exp\{%
{\textstyle\sum_{i=1}^{k}}
\alpha_{i}\log\pi_{i}(\theta\,|\,x)\}$ & when $t=0$\\
$c_{-\infty}(\alpha,m_{\cdot}(x)\pi_{\cdot}(\cdot\,|\,x))\min_{i=1\ldots
,k}m_{i}(x)\pi_{i}(\theta\,|\,x)$ & when $t=-\infty$\\
$c_{\infty}(\alpha,m_{\cdot}(x)\pi_{\cdot}(\cdot\,|\,x))\max_{i=1\ldots
,k}m_{i}(x)\pi_{i}(\theta\,|\,x)$ & when $t=\infty$%
\end{tabular}
\ \ \ \right.
\]
and $m_{t,\alpha}(x)/c_{t}(\alpha,\pi)\leq(\geq)m_{1,\alpha}(x)$ when
$t\leq(\geq)1.$

\noindent Proof: The expressions for $\pi_{t,\alpha}(\cdot\,|\,x)$ for
$t\neq0$ are obvious and%
\begin{align*}
&  \pi_{0,\alpha}(\theta\,|\,x)=c_{0}(\alpha,m_{\cdot}(x)\pi_{\cdot}%
(\cdot\,|\,x))\exp\{%
{\textstyle\sum_{i=1}^{k}}
\alpha_{i}\log m_{i}(x)\pi_{i}(\theta\,|\,x)\}\\
&  =\{\int_{\Theta}\prod\nolimits_{i=1}^{k}m_{i}^{\alpha_{i}}(x)\prod
\nolimits_{i=1i}^{k_{i}}\pi_{i}^{\alpha_{i}}(\theta\,|\,x)\,d\theta
\}^{-1}\prod\nolimits_{i=1}^{k}m_{i}^{\alpha_{i}}(x)\pi_{i}^{\alpha_{i}%
}(\theta\,|\,x)
\end{align*}
so the factor $\prod\nolimits_{i=1}^{k}m_{i}^{\alpha_{i}}(x)$ cancels giving
the result. Finally,
\[
m_{t,\alpha}(x)=c_{t}(\alpha,\pi)\int_{\Theta}\{\sum_{i=1}^{k}\alpha_{i}%
m_{i}^{t}(x)\pi_{i}^{t}(\theta\,|\,x)\}^{1/t}\,d\theta
\]
and this is bounded above (below) by $c_{t}(\alpha,\pi)\int_{\Theta}%
\{\sum_{i=1}^{k}\alpha_{i}m_{i}(x)\pi_{i}(\theta\,|\,x)\}\,d\theta
=c_{t}(\alpha,\pi)m_{1,\alpha}(x)$ when $t\leq(\geq)1$ which gives the
inequality. $\blacksquare$\smallskip

\noindent So the posterior is always proportional to a power mean of the
individual posteriors of the same degree as the power mean of the priors but,
excepting the $t=0$ case, the weights have changed and when $t=-\infty$ or
$t=-\infty$ the prior and\ posterior do not depend on $\alpha.$ The posterior
resulting when $t=1$ is
\begin{equation}
\pi_{1,\alpha}(\theta\,|\,x)=\sum_{i=1}^{k}\left(  \frac{\alpha_{i}m_{i}%
(x)}{\sum_{i=1}^{k}\alpha_{i}m_{i}(x)}\right)  \pi_{i}(\theta\,|\,x)=\sum
_{i=1}^{k}\frac{\alpha_{i}m_{i}(x)}{m_{1,\alpha}(x)}\pi_{i}(\theta\,|\,x),
\label{postlin}%
\end{equation}
and so is a linear combination of the individual posteriors but with different
weights than the prior. The case $t=1$ is called the \textit{linear opinion
pool,} see Stone (1961), and when $t=0$ it is called the \textit{logarithmic
opinion pool}.

The weights staying constant from a priori to a posteriori property for
$\pi_{0,\alpha},$ or even independence from the weights, may seem like an
appealing property but, as discussed in Section 3, these combination rules
have properties that make them inappropriate for combining evidence. A
combination rule is said to be \textit{externally Bayesian} when the rule for
combining the posteriors is the same as the rule for combining the priors. As
shown in Genest (1984a,b), logarithmic pooling is characterized by being
externally Bayesian while linear pooling only satisfies this when there is a
dictatorship, namely, $\alpha_{i}=1$ for some $i,$ as otherwise the weights
differ. Proposition 2 (iii) shows, however, that there is a sense in which
linear pooling can be considered as Bayesian.

Linear pooling has a number of appealing properties.\smallskip

\noindent\textbf{Proposition} \textbf{2.} For Context I, linear pooling
satisfies the following:

\noindent(i) the prior probability measure satisfies the same combination rule
as the densities, namely, $\Pi_{1,\alpha}=\sum_{i=1}^{k}\alpha_{i}\Pi_{i}$ and
similarly for the posterior,

\noindent(ii) marginal priors obtained from $\Pi_{1,\alpha}$ are equal to the
same combination of the marginal priors obtained from the $\Pi_{i}$, and this
is effectively the only rule with this property among all possible combination rules,

\noindent(iii) if $(i,\theta,x)$ is given joint prior distribution with
density $\alpha_{i}\pi_{i}(\theta)f_{\theta}(x)$, then the posterior density
of $\theta$ is given by (\ref{postlin}) and the weight $\alpha_{i}%
m_{i}(x)/m_{1,\alpha}(x)$ is the posterior probability of the index $i$.

\noindent Proof: The proof of (i) is obvious while (ii)\ is proved in McConway
(1981) and holds here with no further conditions. For (iii), note that
$\pi_{i}$ is the conditional prior of $\theta$ given $i$ and $f_{\theta}$ is
the conditional density of $x$ given $\theta.$ Once $x$ is observed the
posterior of $(i,\theta)$ is then given by $\alpha_{i}\pi_{i}(\theta
)f_{\theta}(x)/m_{1,\alpha}(x)$ which implies that the marginal posterior of
$\theta$ is (\ref{postlin}) and the posterior probability of $i$ is
$\alpha_{i}m_{i}(x)/m_{1,\alpha}(x)$. $\blacksquare$\smallskip

\noindent The significance of (i) is that the other combination rules
considered here do not exhibit such simplicity and require more computation to
get the measures. Property (ii) implies that integrating out nuisance
parameters before or after combining does not affect inferences about a
marginal parameter $\psi$ in Context I as the marginal models for $\psi$ are
all the same. Genest (1984c) proves a similar result allowing for negative
$\alpha_{i}$. Property (iii) shows that both the prior $\pi_{1,\alpha}$ and
the posterior $\pi_{1,\alpha}(\cdot\,|\,x)$ can be seen to arise via valid
probability calculations when $\alpha$ is known. A possible interpretation of
this is that $\alpha_{i}$ represents the combiner's prior belief in how well
the $i$-th prior represents appropriate beliefs concerning the true value of
$\theta$ relative to the other priors. The posterior weight $\alpha_{i}%
m_{i}(x)/m_{1,\alpha}(x)$ is then the appropriate modified belief after seeing
the data, as the factor $m_{i}(x)/m_{1,\alpha}(x)$ reflects how well the
$i$-th inference has done at predicting the observed data relative to the
other inference bases. This is a somewhat different interpretation than that
taken by Bunn\ (1981) where $\alpha_{i}$ represents the combiner's prior
belief that the $i$-th inference base is the true one which, in this context,
doesn't really apply.

One commonly cited negative property of linear pooling, see Ladagga (1977), is
that if $A$ and $C$ are independent events for each $\Pi_{i},$ then generally
$\Pi_{1,\alpha}(A\cap C)\neq\Pi_{1,\alpha}(A)\Pi_{1,\alpha}(C).$ It is also to
be noted that if also one of $\Pi_{i}(A)$ or $\Pi_{i\alpha}(C)$ is constant in
$i,$ then independence is preserved and this will be seen to play a role in
linear pooling behaving appropriately when considering statistical evidence,
see Proposition 4 (ii) and the discussion thereafter

\section{Combining Measures of Evidence with Given Prior Weights}

The criterion for choosing an appropriate combination should depend on how
statistical evidence is characterized, as using the evidence to determine
inferences is the ultimate purpose of a statistical analysis. The underlying
idea concerning evidence used here is \textit{the principle of evidence}%
\ which says that there is evidence in favor of the value $\psi$ if $\pi
_{\Psi}(\psi\,|\,x)>\pi_{\Psi}(\psi),$ there is evidence against the value
$\psi$ if $\pi_{\Psi}(\psi\,|\,x)<\pi_{\Psi}(\psi),$ and no evidence either
way if $\pi_{\Psi}(\psi\,|\,x)=\pi_{\Psi}(\psi).$ So, if the data has caused
belief in the value $\psi$ being true to go up, then there is evidence in
favor of this value, etc. The principle of evidence does not require that a
specific numerical measure of evidence be chosen only that any measure used be
consistent with this principle, namely, that there is a cut-off such that the
numerical value greater than (less than) the cut-off corresponds to evidence
in favor of (against) as indicated by the principle. The relative belief ratio
$RB_{\Psi}(\psi\,|\,x)=\pi_{\Psi}(\psi\,|\,x)/\pi_{\Psi}(\psi),$ the ratio of
the posterior to the prior, with the necessary cut-off 1, is used here as it
has a number of good properties, see Evans (2015), and it is particularly
appropriate for the combination of the evidence as easily interpretable
formulas result.

The next result examines the behavior of the combination rules of Section 2
with respect to evidence and is stated initially for the full model parameter
$\theta$.\smallskip

\noindent\textbf{Proposition} \textbf{3}. For Context I, the relative belief
ratio for $\theta$ based on the $\pi_{t,\alpha}$ prior equals
\begin{equation}
RB_{t,\alpha}(\theta\,|\,x)=\frac{m_{1,\alpha}(x)}{m_{t,\alpha}(x)}\sum
_{i=1}^{k}\left(  \frac{\alpha_{i}m_{i}(x)}{\sum_{i=1}^{k}\alpha_{i}m_{i}%
(x)}\right)  RB_{i}(\theta\,|\,x)=\frac{m_{1,\alpha}(x)}{m_{t,\alpha}%
(x)}RB_{1,\alpha}(\theta\,|\,x). \label{rbgen}%
\end{equation}

\noindent Proof: Using $RB_{i}(\theta\,|\,x)=f_{\theta}(x)/m_{i}(x)$ and
$f_{\theta}(x)=\sum_{i=1}^{k}\alpha_{i}f_{\theta}(x),$ then%
\[
RB_{t,\alpha}(\theta\,|\,x)=\frac{f_{\theta}(x)}{m_{t,\alpha}(x)}%
=\frac{m_{1,\alpha}(x)}{m_{t,\alpha}(x)}\sum_{i=1}^{k}\frac{\alpha_{i}%
m_{i}(x)}{m_{1,\alpha}(x)}RB_{i}(\theta\,|\,x).\text{ }\blacksquare
\]

\noindent This result shows the value of using the relative belief ratio to
express evidence since the combination rule, at least for power means, is
quite simple and natural. Notice too that if there are only $l$ distinct
priors, then the combination rules for the priors, posteriors and relative
belief ratios are really only based on these distinct priors and the weights
change only by summing the $\alpha_{i}$ that correspond to common priors.

The result in Proposition 3 is another indication that the correct way to
combine priors, from the point of view of measuring evidence, is via linear
pooling as $RB_{t,\alpha}(\theta\,|\,x)$ is always proportional to
$RB_{1,\alpha}(\theta\,|\,x).$ The constant multiplying $RB_{1,\alpha}%
(\theta\,|\,x)$ in (\ref{rbgen}) suggests that finding $t$ that minimizes
$m_{1,\alpha}(x)/m_{t,\alpha}(x)$, leads to the power mean prior that
maximizes the amount of mass the prior places at $\theta_{true},$ see
Proposition 7 (iv) . But there is a significant reason for preferring
$RB_{1,\alpha}(\theta\,|\,x)$ over the other possibilities. For suppose that
$RB_{i}(\theta\,|\,x)<1$ for all $i\ $or\ $RB_{i}(\theta\,|\,x)>1$ for all
$i.$ Then it is clear that $RB_{1,\alpha}(\theta\,|\,x)<1$ in the first case
and $RB_{1,\alpha}(\theta\,|\,x)>1$ in the second case. In the first case
there is a consensus that there is evidence against $\theta$ being the true
value and in the second case there is a consensus that there is evidence in
favor of $\theta$ being the true value. In other words $RB_{1,\alpha}$ is
consensus preserving and this seems like a necessary property for any approach
to combining evidence.

A formal definition is now provided which takes into account that sometimes
$RB_{i}(\theta\,|\,x)=1$ indicating that there is no evidence either way and
the $i$-th inference base is agnostic about whether or not $\theta$ is the
true value.\smallskip

\noindent\textbf{Definition} A rule for combining evidence about a parameter
is called \textit{consensus preserving }if, whenever at least one of the
inference bases indicates evidence in favor of (against) a value of the
parameter and the remaining inference bases do not give evidence against (in
favor), then the rule gives evidence in favor of (against) the value and if no
inference base indicates evidence one way or the other then neither does the
combination.\smallskip

\noindent The following property is immediately obtained for linear
pooling.\smallskip

\noindent\textbf{Proposition} \textbf{4}. For Context I, whenever $\alpha
_{i}>0$ for all $i,$ then (i) $RB_{1,\alpha}$ is consensus preserving and (ii)
whenever $RB_{i}(\theta\,|\,x)\leq(\geq)1$ for all $i,$ then $RB_{1,\alpha
}(\theta\,|\,x)=1$ iff $RB_{i}(\theta\,|\,x)=1$ for all $i$.\smallskip

\noindent The property of preserving consensus is similar to the unanimity
principle for priors, see Clemen and Winkler (1999), which says that if all
the priors are the same, then the combination rule must give back that prior.

Proposition 4 (ii) indicates that linear pooling deals correctly with
independent events at least with respect to evidence. For note that, for
probability measure $P$ and events $A$ and $C$ satisfying $P(A\cap C)>0,$ then
$A$ and $C$ are statistically independent iff $RB(A\,|\,C)=P(A\cap
C)/P(A)P(C)=1$ and independence is equivalent to saying that the occurrence of
$C$ provides no evidence concerning the truth or falsity of $A.$ Now consider
the statistical context and suppose $RB_{i}(\theta\,|\,x)=1$ and further
suppose that all the probabilities are discrete. This implies that $f_{\theta
}(x)=m_{i}(x)$ which implies the joint prior density at $(\theta,x)$ factors
as $f_{\theta}(x)\pi_{i}(\theta)=m_{i}(x)\pi_{i}(\theta)$ and so the events
$\{\theta\}$ and $\{x\}$ are statistically independent in the $i$-th inference
base. If this holds for each $i,$ then $m_{i}(x)$ is constant in $i$ and so
indeed $RB_{1,\alpha}(\theta\,|\,x)=1$ implies that these events are
independent when the prior is the linear pool. With a continuous prior then
$RB_{i}(\theta\,|\,x)=1$ can also happen, but typically this event has prior
probability 0.

It is of interest to determine whether or not any of the other rules based on
the means are consensus preserving. The inequality in Proposition 1 and
Proposition 3 imply that, when $t\leq1,$ then $RB_{t,\alpha}(\theta\,|\,x)\geq
RB_{1,\alpha}(\theta\,|\,x)/c_{t}(\alpha,\pi)$ and since $c_{t}(\alpha
,\pi)\geq1,$ with the inequality typically strict when $t<1,$ this suggests
that $RB_{t,\alpha}$ might even contradict the consensus of evidence in favor.
A similar argument holds for $t>1.$ The following example shows that generally
the combination rules based on power means of priors aren't consensus
preserving.\smallskip

\noindent\textbf{Example 1. }\textit{Power means of priors aren't generally
consensus preserving.}

Suppose $\mathcal{X}=\{0,1\},\Theta=\{a,b\},f_{a}(0)=1/4,f_{b}(0)=1/3$ and
$x=0$ is observed. There are two priors given by $\pi_{1}(a)=p_{1}$ and
$\pi_{2}(a)=p_{2}.$ Then $m_{1}(0)=(4-p_{1})/12,m_{2}(0)=(4-p_{2}%
)/12,RB_{1}(a\,|\,0)=3/(4-p_{1}),RB_{2}(a\,|\,0)=3/(4-p_{2})$ so both
inference bases give evidence against when $p_{1}<1,p_{2}<1$ and
$RB_{i}(a\,|\,0)=1$ when $p_{i}=1$ so no evidence either way is obtained from
the data when a statistician is categorical in their beliefs. Note being
categorical in your beliefs is a possible choice provided it doesn't lead to
prior-data conflict. When $\alpha=1/2,$ so the two priors are being given
equal weight, then $\pi_{1,1/2}(a)=(p_{1}+p_{2})/2,m_{1,1/2}(0)=(m_{1}%
(0)+m_{2}(0))/2=(8-p_{1}-p_{2})/24$ and%
\[
RB_{1,1/2}(a\,|\,0)=(m_{1}(0)RB_{1}(a\,|\,0)+m_{2}(0)RB_{2}%
(a\,|\,0))/2m_{1,1/2}(0)=6/(8-p_{1}-p_{2}).
\]
When $p_{1}=1/4,p_{2}=1,$ so statistician 2 is categorical in their beliefs,
then $m_{1,1/2}(0)=9/32,RB_{1}(a\,|\,0)=4/5,RB_{2}(a\,|\,0)=1$ and
$RB_{1,1/2}(a\,|\,0)=8/9$ so there is evidence against.

Now consider $\pi_{0,\alpha}$ given by $\pi_{0,\alpha}(a)=p_{1}^{\alpha}%
p_{2}^{\alpha}/(p_{1}^{\alpha}p_{2}^{\alpha}+(1-p_{1})^{1-\alpha}%
(1-p_{2})^{1-\alpha}).$ So $m_{0,\alpha}(0)=\pi_{0,\alpha}(a)/4+\pi_{0,\alpha
}(b)/3$ and when $p_{1}=1/4,p_{2}=1,$ then $m_{0,1/2}(0)=1/4$ and so
\[
RB_{0,1/2}(a\,|\,0)=\frac{m_{1,1/2}(0)}{m_{0,1/2}(0)}RB_{1,1/2}(a\,|\,0)=\frac
{9/32}{1/4}8/9=1
\]
and so $RB_{0,1/2}(a\,|\,0)$ indicates no evidence either way. Therefore, the
case $t=0$ is not consensus preserving.

Next consider the case $t=-\infty,$ so $\pi_{-\infty,\alpha}(a)=\min
\{p_{1},p_{2}\}/(\min\{p_{1},p_{2}\}+\min\{1-p_{1},1-p_{2}\}).$ When
$p_{1}=1/4,$ then $\pi_{-\infty,\alpha}(a)=1$ and so $m_{-\infty,\alpha
}(0)=1/4$ which implies that
\[
RB_{-\infty,1/2}(a\,|\,0)=\frac{m_{1,1/2}(0)}{m_{-\infty,1/2}(0)}%
RB_{1,1/2}(a\,|\,0)=\frac{9/32}{1/4}8/9=1.
\]
So also the case $t=-\infty$ is not consensus preserving.

If there is evidence against (in favor of) an event, then a property of the
relative belief ratio gives that there is evidence in favor of (against) its
complement and, if there is no evidence either way for an event, then there is
no evidence either way for its complement, see Evans (2015), Proposition 4.2.3
(i). So in this example the priors $\pi_{0,\alpha}$ and $\pi_{-\infty,\alpha}$
also do not preserve consensus with respect to $\theta=b.$ $\blacksquare
$\smallskip

So far no case has been found where a combination based on a power mean
actually reverses a consensus and it is a reasonable conjecture, based on many
examples, that this will never happen but a proof is not obvious. It is still
disturbing, however, that, as Example 1 illustrates, the inference bases can
be equally weighted with none giving evidence in favor and at least one giving
evidence against but the determination is that no evidence is obtained either
way. It can be argued that a rule that behaves like this is not reflecting
what the overall conclusion should be about the evidence.

There is another interesting consequence of Proposition 3 which is relevant
when the goal is to estimate $\theta.$ The natural estimate is the relative
belief estimate $\theta(x)=\arg\sup_{\theta}RB(\theta\,|\,x)$ and the accuracy
of this estimate is assessed by the \textit{plausible region} $Pl(x)=\{\theta
:RB(\theta\,|\,x)>1\},$ the set of values for which there is evidence in
favor. For example, the "size" of $Pl(x)$ and its posterior content together
provide an a posteriori indication of how accurate $\theta(x)$ is. Ideally we
want $Pl(x)$ "small" and its posterior content high. Note that it is easy to
show in general that $RB(\theta(x)\,|\,x)>1$ so $\theta(x)\in Pl(x)$ provided
$RB(\theta\,|\,x)$ is not 1 for all $\theta$ which only occurs when the data
indicates nothing about the true value.\smallskip

\noindent\textbf{Corollary} \textbf{5}. Whenever $RB_{t,\alpha}(\theta\,|\,x)$
is not 1 for all $\theta$ and $\alpha_{i}>0$ for all $i,$ then $\arg
\sup_{\theta}RB_{t,\alpha}(\theta\,|\,x)=\arg\sup_{\theta}RB_{1,\alpha}%
(\theta\,|\,x).$\smallskip

\noindent So the estimate of $\theta$ based on maximizing the evidence is
always determined by linear pooling. It is not the case, however, that the
plausible region is independent of $t$ because of the constant $m_{1,\alpha
}(x)/m_{t,\alpha}(x).$

The following underscores the role of linear pooling in preserving
consensus.\smallskip

\noindent\textbf{Corollary} \textbf{6}. The set $\{\theta:RB_{i}%
(\theta\,|\,x)>1$ for all $i\}=\cap_{i=1}^{k}Pl_{i}(x)\subset Pl_{1,\alpha
}(x)$ and $\min_{i}\Pi_{i}(Pl_{i}(x)\,|\,x)\leq\Pi_{1,\alpha}(Pl_{1,\alpha
}(x)\,|\,x)\leq\max_{i}\Pi_{i}(Pl_{i}(x)\,|\,x).$\smallskip

\noindent So the set of $\theta$ where there is a consensus that there is
evidence in favor is always contained in the plausible region determined by
linear pooling. A similar comment applies to the \textit{implausible region}
which is the set of all values where there is evidence against. While it might
be tempting to quote the region $\cap_{i=1}^{k}Pl_{i}(x),$ there is no
guarantee that any of the relative belief estimates will be in this set,
whether determined by $RB_{1,\alpha}(\cdot\,|\,x)$ or any of the $RB_{i}%
(\cdot\,|\,x).$

The situation with respect to the assessment of the hypothesis $H_{0}%
:\theta=\theta_{0}$ is a bit different. Clearly, if $RB_{i}(\theta
_{0}\,|\,x)>(<)1$ for all $i,$ so there is a consensus that there is evidence
in favor of (against) $H_{0},$ then $RB_{1,\alpha}(\theta_{0}\,|\,x)$
preserves this consensus. In general, if one assesses the evidence for $H_{0}$
via a relative belief ratio $RB(\theta_{0}\,|\,x),$ then the posterior
probability $\Pi(RB(\theta_{0}\,|\,x)\leq RB(\theta_{0}\,|\,x)\,|\,x)$ can be
taken as a measure of the strength of the evidence, see Evans (2015). In the
context under discussion here, it follows from (\ref{rbgen}) that the event
$\{\theta:RB_{t,\alpha}(\theta\,|\,x)\leq RB_{t,\alpha}(\theta_{0}%
\,|\,x)\}=\{\theta:RB_{1,\alpha}(\theta\,|\,x)\leq RB_{1,\alpha}(\theta
_{0}\,|\,x)\}$ for all $t.$ Of course, the posterior probability of this event
will depend on $t$ but linear pooling completely determines the event.

Now suppose that interest is in the quantity $\psi=\Psi(\theta)$ and the
assumptions of Context I hold so that prior beliefs only differ concerning the
value of $\psi$ which implies that the inference bases only differ with
respect to the priors on $\psi.$ This situation may arise when the analysts
all agree to use a common default prior on the nuisance parameters. Then we
can treat $\psi\,$\ as the model parameter for the common model $\{m_{\psi
}:\psi\in\Psi\}$ and the relevant linear pooling rule is%
\begin{equation}
RB_{1,\alpha,\Psi}(\psi\,|\,x)=\sum_{i=1}^{k}\frac{\alpha_{i}m_{i}%
(x)}{m_{1,\alpha}(x)}RB_{i,\Psi}(\psi\,|\,x) \label{linrule}%
\end{equation}
and all the results derived for $\theta$ apply.

In general it can be expected that some inference bases will indicate evidence
in favor of $\psi$ being the true value and some will indicate evidence
against, but $RB_{1,\alpha,\Psi}(\psi\,|\,x)$ will indicate evidence one way
or the other. This depends on the values assumed by the $RB_{i,\Psi}%
(\psi\,|\,x)$ as well as the weights $\alpha_{i}m_{i}(x)/m_{1,\alpha}(x).$ It
is to be noted, however, that there is nothing paradoxical about this as it is
possible that there is evidence in favor of a set $A$ and contrary evidence
for a set $B$ even though $A\subset B.$ This situation is discussed in
Corollary 4.2.1 of Evans (2015)\ and it is shown there that, when the prior
probabilities are taken into account, the apparent paradox disappears. For
example, if the prior probability of $A$ is a small relative to the prior
probability of $B,$ then such a reversal can occur and so the evidence
measures disagreeing is not paradoxical. As an example of this, if there is
evidence that a small subgroup of a population has extreme views on an issue,
this is not evidence favoring the whole population having such extreme views
and, on balance, the evidence could well indicate otherwise depending on the
relative size of the subgroup. In essence, measuring evidence is quite
different than measuring beliefs via probability as it is the change in belief
from a priori to a posteriori that informs us about the evidence the data is expressing.

Consider now the context where $\breve{x}_{n}=(x_{1},\ldots,x_{n})$ is an
$i.i.d.$ sample. The following result gives the consistency of this approach
when the basic sample space for the $x_{i}$ and the model parameter space
$\Theta$ are finite. Such results will hold more generally but require some
mathematical constraints on densities and this is not pursued further here.
Let $\psi_{1,\alpha}(\breve{x}_{n})=\arg\max_{\psi}RB_{\Psi,1,\alpha}%
(\psi\,|\,\breve{x}_{n})$ be the relative belief estimate of $\psi$ based on
linear pooling. All the convergence results are almost everywhere as
$n\rightarrow\infty$ with the proof in the Appendix$.$\smallskip

\noindent\textbf{Proposition} \textbf{7}. For Context I, suppose $\breve
{x}_{n}=(x_{1},\ldots,x_{n})$ is an $i.i.d.$ sample from a distribution in a
model having a finite parameter space $\Theta$ and each prior for $\theta$\ is
everywhere positive on $\Theta.$ Then \newline\noindent(i) $RB_{1,\alpha,\Psi
}(\psi_{0}\,|\,\breve{x}_{n})\rightarrow I_{\{\psi_{true}\}}(\psi_{0}%
)/\pi_{1,\alpha,\Psi}(\psi_{0})$ and $\Pi_{1,\alpha,\Psi}(RB_{1,\alpha,\Psi
}(\psi\,|\,\breve{x}_{n})\leq RB_{1,\alpha,\Psi}(\psi_{0}\,|\,\breve{x}%
_{n}))\,|\,\breve{x}_{n})\rightarrow I_{\{\psi_{true}\}}(\psi_{0}),$%
\newline\noindent(ii) $\psi_{1,\alpha}(\breve{x}_{n})\rightarrow\psi
_{true},Pl_{1,\alpha,\Psi}(\breve{x}_{n})\rightarrow\{\psi_{true}\}$ and
$\Pi_{1,\alpha,\Psi}(Pl_{1,\alpha,\Psi}(\breve{x}_{n})\,|\,\breve{x}%
_{n})\rightarrow1,$

\noindent(iii) $\alpha_{i}m_{i}(x)/m_{1,\alpha}(x)\rightarrow\alpha_{i}\pi
_{i}(\theta_{true})/\pi_{1,\alpha}(\theta_{true}),$

\noindent(iv) $m_{1,\alpha}(x)/m_{t,\alpha}(x)\rightarrow\pi_{1,\alpha}%
(\theta_{true})/\pi_{t,\alpha}(\theta_{true}).$\smallskip

\noindent Noting that when $1/\pi_{1,\alpha,\Psi}(\psi_{0})>1,$ then
Proposition 7 (i) says that the evidence in favor of (against) $H_{0}%
:\Psi(\theta)=\psi_{0},$ based on the combination, goes to categorical when
$H_{0}$ is true (false). Part (ii) says that the relative belief estimate
based on the combination is consistent. Part (iii) implies that, when the
priors are equally weighted, then the inference base whose prior gives the
largest value to the true value will inevitably have the largest weight in
determining the combined evidence.

Overall, the conclusion reached here is that linear pooling is the most
natural way to combine evidence at least among the power means. As such,
attention is restricted to this case hereafter. Various authors, when
discussing the combination of priors, have come to a similar conclusion. For
example, O'Hagan et al. (2006), when considering the full spectrum of methods
for combining priors, write the following, "In general, it seems that a
simple, equally weighted, linear opinion pool is hard to beat in practice."
The results developed here support such a conclusion when considering evidence.

\section{Determining the Prior Weights}

The discussion so far has assumed that $\alpha$ is known. So arguments or
methodologies for choosing $\alpha$ need to be considered. There are several
possible approaches to determining a suitable choice of the prior weights and
nothing novel is proposed here. As previously mentioned, the $\alpha_{i}$ can
represent the combiner's beliefs concerning how well the $i$-th prior
represents appropriate beliefs about $\theta.$ The combiner's beliefs should
of course be based upon experience or knowledge concerning the various
proposers of the priors. In absence of such knowledge then uniform weights,
namely, $\alpha=(1/k,\ldots,1/k),$ seem reasonable. Genest and McConway (1990)
provides a good survey of various approaches to choosing $\alpha$. Also,
DeGroot (1974) and Lehrer and Wagner (1981) present a novel iterative approach
to determining a consensus $\alpha$ among the proposers.

In Context I notice that the weights $\alpha_{i}m_{i}(x)/m_{1,\alpha}(x)$ only
depend on the data through some function of the value of the minimal
sufficient statistic (mss) for the model. So, for example, if the priors are
distinct and equally weighted via $\alpha=(1/k,\ldots,1/k),$ then the weight
of the $i$-th prior is $m_{i}(x)/(m_{1}(x)+\cdots+m_{k}(x))$ and so more
weight is given to those inference bases that do a better job, relatively
speaking, of predicting a priori the observed value of this function of the
mss). Since it is only the observed value of the mss that is relevant for
inference, this seems sensible. There is the possibility, however, to weight
some priors more than others for a variety of reasons.

A prior can also be placed on $\alpha,$ the results examined for a number
different choices of $\alpha$ and summarized in a way that addresses the issue
of whether or not the inferences are sensitive to $\alpha.$ For example,
suppose the goal is to determine if there is evidence for or against the
hypothesis $H_{0}:\Psi(\theta)=\psi_{0}.$ For a given weighting $\alpha_{0}$
the evidence for or against will be determined by the value $RB_{1,\alpha
_{0},\Psi}(\psi_{0}\,|\,x).$ Accordingly, a Dirichlet prior with mode at
$\alpha_{0}$ and with some degree of concentration around this point could be
used to assess the robustness of the combination inferences. In particular,
for each generated value of $\alpha$ from the prior, one can record whether
evidence in favor of or against $H_{0}$ was obtained together with the
strength of the evidence. If a great proportion of the results gave results
similar to those obtained with the weights $\alpha_{0},$ then this would
provide some assurance that the conclusions drawn are robust to deviations. A
similar approach can be taken to estimation problems where the relative belief
estimate is given by $\psi(x)=\arg\sup_{\psi}RB_{1,\alpha_{0},\Psi}%
(\psi\,|\,x).$ When $\Psi$ is 1-dimensional then a histogram of the estimates
obtained in the simulation and histograms of the prior and posterior contents
of $Pl_{\Psi}(x)$ will provide an indication of the dependence on $\alpha
_{0}.$

\section{The General Problem}

The general Context II is more complicated and an overall solution is not
proposed here rather a special case is considered when there is a common data
set $x.$ Since Context II covers Context I, it is necessary that any rule
proposed for such situations agrees with what is determined for Context I when
that applies.

While the prior on $\psi$ can be taken to be the mixture $\sum_{i=1}^{k}%
\alpha_{i}\pi_{i,\Psi_{i}},$ the overall posterior does not have a clear
definition as it is not obvious how to form the likelihood. For example, the
model parameters $\theta_{i}$ may not be comparable even if the parameter of
interest $\psi=\Psi(\theta_{i})$ always references the same real-world object
as is assumed here. In some contexts there may be good arguments for a
specific definition of the likelihood but this issue is not addressed further
as our focus is on combining the evidence as characterized by the individual
inference bases.

The simplest approach to characterizing the evidence in Context II is to use
the obvious expression%
\begin{equation}
RB_{1,\alpha,\Psi}^{\ast}(\psi\,|\,x)=\sum_{i=1}^{k}\frac{\alpha_{i}m_{i}%
(x)}{m_{1,\alpha}(x)}RB_{i,\Psi}(\psi\,|\,x), \label{linrule2}%
\end{equation}
where again $RB_{i,\Psi}(\psi\,|\,x)$ and $m_{i}(x)$ arise from the $i$-th
inference base and $m_{1,\alpha}(x)=\sum_{i=1}^{k}\alpha_{i}m_{i}(x_{i})$.
This will agree with the answer obtained in Context I when it applies, but
generally $RB_{1,\alpha,\Psi}^{\ast}(\psi\,|\,x)$ is not the ratio of the
posterior of $\psi$ to its prior and as such it cannot be claimed that it is a
valid characterization of the evidence as $RB_{1,\alpha,\Psi}(\psi\,|\,x)$ is
in\ Context I.

One approach to defining a posterior in Context II is to use the argument
known as Jeffrey conditionalization, see Diaconis and Zabell (1982). This
involves considering the probabilities on the partition given by
$i\in\{1,\ldots,k\}$ completely separately from the probabilities on $\Psi$
given $i.$ In this scenario the probabilities on the partition elements $i$
are updated differently than the probabilities given a partition element. This
makes some sense here because one can think of the $\alpha_{i}$ as being the
combiner's prior probabilities concerning the relevance of the $i$-th
inference base to inference about $\psi$ and these are separate from the
individual analyst's beliefs expressed about the true value of $\psi.$ For
example, some of the inference bases could be formed by teams with much more
relevant experience than some of the others and so be more heavily weighted.
As in Proposition 2 (iii), $\alpha_{i}m_{i}(x)$ can be thought of as the prior
probability of $(i,x)$ so, after observing $x,$ the posterior probability of
$i$ is again given by $\alpha_{i}m_{i}(x)/m_{1,\alpha}(x)$ and given $(i,x)$
the posterior of $\psi$ is $\pi_{i,\Psi}(\psi\,|\,x_{i}).$ Both of these arise
via regular (Bayesian) conditionalization but with priors on different
objects. Using the Jeffrey's conditionalization argument, the overall
posterior of $\psi$ is%
\begin{equation}
\pi_{1,\alpha,\Psi}^{\ast}(\psi\,|\,x)=\sum_{i=1}^{k}\frac{\alpha_{i}m_{i}%
(x)}{m_{1,\alpha}(x)}\pi_{i,\Psi}(\psi\,|\,x). \label{Jeffreypost}%
\end{equation}
It is not the case, however, that in general (\ref{linrule2}) results as the
relative belief ratio formed from (\ref{Jeffreypost}) and the prior
$\pi_{1,\alpha,\Psi}$ although (\ref{Jeffreypost}) will be used to determine
probabilities like the content of the plausible region determined by
(\ref{linrule2}). One reason for not using $\pi_{1,\alpha,\Psi}$ and
$\pi_{1,\alpha,\Psi}^{\ast}(\cdot\,|\,x)$ to determine the evidence is that
this does not lead to a linear pooling of the characterizations of the
evidence via the individual relative belief ratios and so the nice properties
of linear pooling are lost, see Example 3. Notice that (\ref{linrule2}) will
satisfy all the properties of linear pooling established for Context I with
the exception of Proposition 2 (iii) and Jeffreys conditionalization is then
used to justify the mixture. In particular, $RB_{1,\alpha,\Psi}^{\ast}%
(\psi\,|\,x)$ will preserve a consensus about evidence in favor or against.

The following result characterizes what happens as sample size grows and is
proved in the Appendix.\smallskip

\noindent\textbf{Proposition} \textbf{8}. Suppose $\breve{x}_{n}=(x_{1}%
,\ldots,x_{n})$ is an $i.i.d.$ sample from a distribution in at least one of
the models and each of the parameter spaces $\Theta_{i}$ is finite with the
prior everywhere positive on $\Theta_{i}.$ Denoting the set of indices
corresponding to the models containing the true distribution by $J,$ then as
$n\rightarrow\infty$\newline\noindent(i) $\alpha_{i}m_{i}(\breve{x}%
_{n})/m_{1,\alpha}(\breve{x}_{n})\rightarrow w_{i}=I_{J}(i)\alpha_{i}\pi
_{i}(\theta_{itrue})/\sum_{j\in J}\alpha_{j}\pi_{j}(\theta_{jtrue})\geq0$ and
$\sum_{i=1}^{k}w_{i}=1\ ,$ \newline(ii) $RB_{1,\alpha,\Psi}^{\ast}%
(\psi\,|\,\breve{x}_{n})\rightarrow I_{\{\psi_{true}\}}(\psi)\sum_{i=}%
^{k}w_{i}/\pi_{i\Psi}(\psi)$ which is greater than $1$ when $\psi=\psi_{true}$
\newline(iii) $\lim\pi_{1,\alpha,\Psi}^{\ast}(\psi\,|\,\breve{x}%
_{n})=I_{\{\psi_{true}\}}(\psi).$\smallskip

\noindent So Proposition 8 shows that $RB_{1,\alpha,\Psi}^{\ast}%
(\cdot\,|\,\breve{x}_{n})$ and $\pi_{1,\alpha,\Psi}^{\ast}(\cdot\,|\,\breve
{x}_{n})$ provide consistent inferences and the weights converge to
appropriate values.

There is another significant difference between (\ref{linrule2}) and
(\ref{linrule}). In Context I the weights all depended on the data through the
same function of a constant mss for the full common model. Furthermore, if
$A(x)$ is an ancillary statistic for the full model, then it is seen that the
$i$-th weight satisfies $\alpha_{i}m_{i}(x)/m_{1,\alpha}(x)=\alpha_{i}%
m_{i}(x\,|\,A(x))/m_{1,\alpha}(x\,|\,A(x)).$ This implies that the weights are
comparable as they are all concerned with predicting essentially the same data
and moreover they are not concerned with predicting aspects of the data that
have no relation to the quantity of interest. In Context II this is not the
case which raises the question of whether or not the weights are comparable.

It is not obvious how to deal with this issue in general, but in some contexts
the structure of the models is such that $x\leftrightarrow(L(x),A(x))$ where
$L$ has fixed dimension and $A$ is ancillary for each model. For example, if
all the models are location models, then $x=(x_{1},\ldots,x_{n})^{\prime}%
=\bar{x}1_{n}+A(x),$ where $1_{n}$ is a column of 1's, and $A(x)=(x_{1}%
-\bar{x},\ldots,x_{n}-\bar{x})^{\prime}$ is ancillary. In such a case it is
desirable to determine the weights based on how well the inference bases
predict the value of $L(x)$ and not $A(x).$ To take account of this it is
necessary that Jeffrey conditionalization be modified so that the $i$-th
weight is now proportional to $\alpha_{i}m_{i}(x\,|\,A(x))$ where $m_{i}%
(\cdot\,|\,A(x))$ is the $i$-th prior predictive of the data given $A(x).$
Examples 4 and 5 illustrate this modification.

While Proposition 8 does not apply with the conditional weights, a similar
result can be proved\ and for this some assumptions are imposed to simplify
the proof. Let the basic sample space be such that there is a finite ancillary
partition $(B_{1},\ldots,B_{m}),$ applicable to each of the $k$ models, and
for any $n$ the ancillary is given by $A(\breve{x}_{n})=(n_{1}(\breve{x}%
_{n}),\ldots,n_{m}(\breve{x}_{n}))$ where $n_{i}(\breve{x}_{n})$ records the
number of values in the sample that lie in $B_{i}.$ Then the probability
distribution of $A(\breve{x}_{n})$ for the $i$-th model is given by the
multinomial$(n,p_{i1},\ldots,p_{im})$ where the $p_{ij}$ are fixed and
independent of the model parameter and denote this probability function at the
observed data by $f_{i}(\breve{n}(\breve{x}_{n}))$ where $\breve{n}(\breve
{x}_{n})=(n_{1}(\breve{x}_{n}),\ldots,n_{m}(\breve{x}_{n})).$ Suppose that
each parameter space $\Theta_{i}$ is finite with the prior $\pi_{i}$
everywhere positive. Let $\alpha_{i}\propto\alpha_{i}^{\ast}/f_{i}(\breve
{n}(\breve{x}_{n}))$ for $\alpha^{\ast}\in S_{k}$ and $J$ denote the set of
indices containing the true distribution. Calling these requirements condition
$\bigstar,$ the following is proved in the Appendix.\smallskip

\noindent\textbf{Proposition} \textbf{9}. If condition $\bigstar$ holds, then
$\alpha_{i}m_{i}(\breve{x}_{n})/m_{1,\alpha}(\breve{x}_{n})\rightarrow
w_{i}=I_{J}(i)\alpha_{i}^{\ast}\pi_{i}(\theta_{itrue})/\sum_{j\in J}\alpha
_{j}^{\ast}\pi_{j}(\theta_{jtrue})\geq0$ and $\sum_{i=1}^{k}w_{i}%
=1.$\smallskip

\noindent Proposition 9 provides the desirable consistency result as the only
thing that is affected here are the weights which have been shown to have the
correct asymptotic property.

Of course, this result needs to be generalized to handle even a situation like
the location model. For this some conditions on the models and priors are
undoubtedly required but this is not pursued further here. One key component
of the proof is the existence of the ancillary partition $(B_{1},\ldots
,B_{m})$ and such a structural element seems necessary generally to get
comparability of the weights. In group-based models, like linear regression
and many others, such a structure exists via the usual ancillaries, see
Example 5. As an approximation, a finite ancillary partition can be
constructed via the ancillary statistic in question and so Proposition 9 is
applicable. It should also be noted that, if the original model is replaced by
the conditional model given the ancillary, then (\ref{linrule2}) gives the
same answer as this modification as the values of $RB_{i,\Psi}(\psi\,|\,x)$
are unaffected by the conditioning.

\section{Examples}

Some examples are now considered and initially a very simple context is
considered where understanding of what is going on is enhanced by the
availability of closed form expressions for relevant quantities. \smallskip

\noindent\textbf{Example 2. }\textit{Location-normal model with normal
priors.}

Suppose $x=(x_{1},\ldots,x_{n})$ is a sample from a $N(\mu,\sigma_{0}^{2})$
distribution where the mean is unknown but the variance is known. It might be
more appropriate to model this with an unknown variance but this situation
will suffice for illustrative purposes. The model is then given by, after
reducing to a mss, the collection of $N(\mu,\sigma_{0}^{2}/n)$ distributions
and so this is Context I. Suppose there are three analysts and they express
their priors for $\mu$ as $N(\mu_{i},\tau_{i}^{2})$ distributions for
$i=1,2,3\ $so the posteriors are $N((n/\sigma_{0}^{2}+1/\tau_{i}^{2}%
)^{-1}(n\bar{x}/\sigma_{0}^{2}+\mu_{i}/\tau_{i}^{2}),(n/\sigma_{0}^{2}%
+1/\tau_{i}^{2})^{-1})$ and these determine the relative belief ratios. For
combining, the prior predictives are also needed and the $i$-th prior
predictive density $m_{i}$ for $\bar{x}$ is the $N(\mu_{i},\sigma_{0}%
^{2}/n+\tau_{i}^{2})$ density. Suppose the inference bases are equally
weighted so the posterior weight of the $i$-th analysis relative to the others
is determined by how well the observed value $\bar{x}$ fits the $N(\mu
_{i},\sigma_{0}^{2}/n+\tau_{i}^{2})$ distribution. Note, however, that even if
there is a perfect fit, in the sense that $\bar{x}=\mu_{i},$ the weight still
depends on the quantity $\sigma_{0}^{2}/n+\tau_{i}^{2}.$ For example, if the
$\mu_{i}$ are all equal and there is a perfect fit, then the $i$-th weight is
proportional to $(1+n\tau_{i}^{2}/\sigma_{0}^{2})^{-1/2}$ and this weight goes
to 0 as $\tau_{i}^{2}\rightarrow\infty$ with the other prior variances
constant and goes to its biggest value when $\tau_{i}^{2}\rightarrow0.$ This
suggests that making a prior quite diffuse leads to reducing the impact such a
prior has in the combined analysis.

\begin{figure}[t]%
\centering
\includegraphics[
height=2.0868in,
width=3.474in
]%
{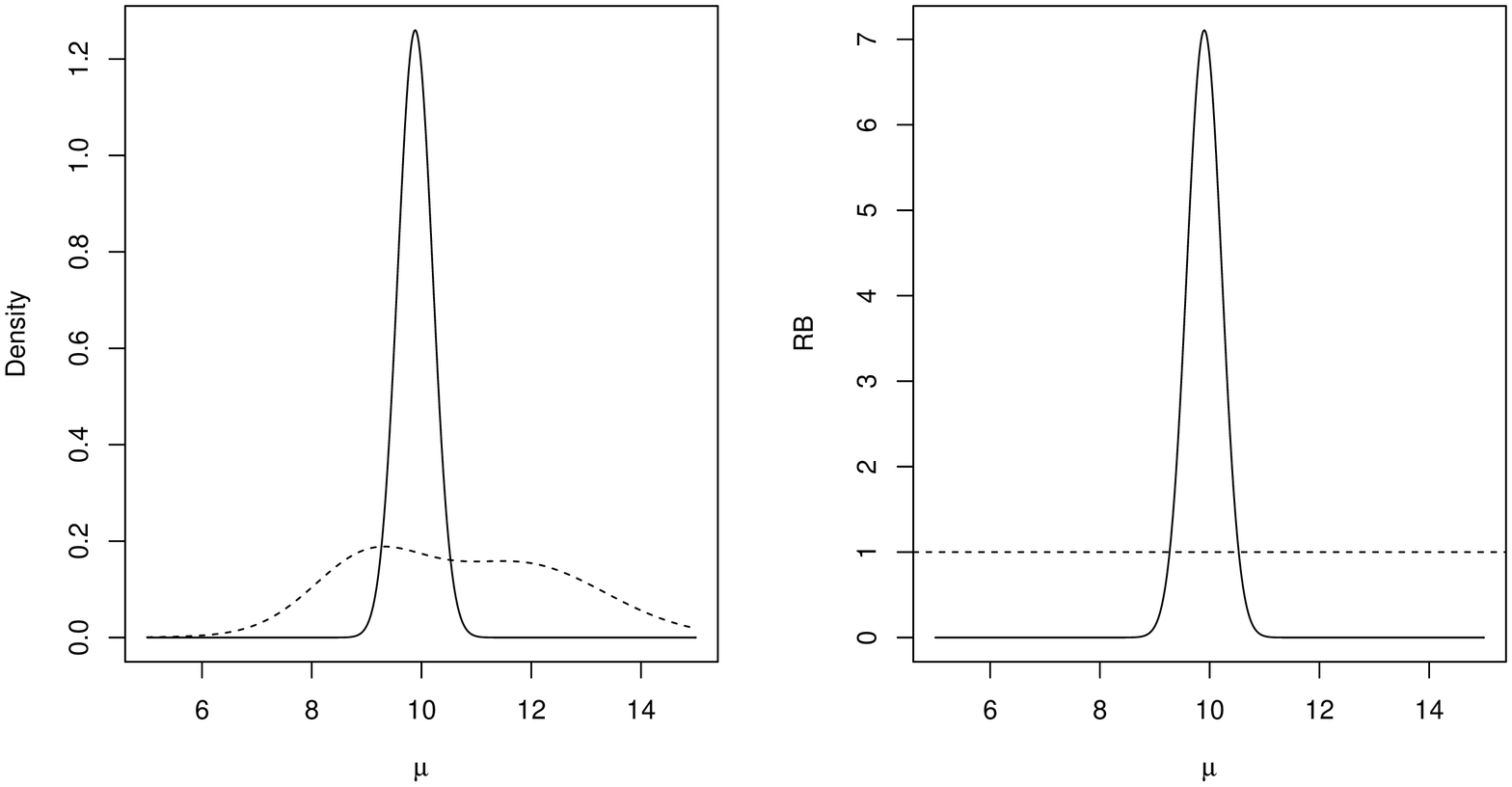}%
\caption{Plots of prior (- - -) and posterior (---) densities in the left
panel and, in the right panel, plots of the relative belief ratio (---) for
$\mu$ and constant 1 (- - -) in Example 2 when $n=10.$}%
\label{fig1}%
\end{figure}

As a specific data example, suppose the true value of $\mu=10,$ with
$\sigma_{0}=1$ and consider the results for several sample sizes
$n=5,10,25,100.$ Data were generated from the true distribution obtaining the
values $\bar{x}=10.92,9.87,9.96,10.12$ respectively. For the priors consider
$(\mu_{1},\tau_{1}^{2})=(12,2),(\mu_{2},\tau_{2}^{2})=(9,1),(\mu_{3},\tau
_{3}^{2})=(11,4)$ with the priors weighted equally. Figure \ref{fig1} plots
the combined prior, posterior and relative belief ratio for the $n=10$ case.
Table \ref{tab1} records the estimates of $\mu,$ the plausible regions
together with the posterior and prior contents of these intervals for each
inference base and linear pooling. Note that in this case, because the model
is the same for each inference base and $\mu$ is the model parameter, the
estimates are all equal to the MLE of $\mu$ but the plausible intervals and
their posterior contents differ. $\blacksquare$\smallskip%
%

\begin{table}[tp] \centering
\begin{tabular}
[c]{|c|c|c|c|c|c|}\hline
$n$ & est. & $\mathcal{I}_{1}$ & $\mathcal{I}_{2}$ & $\mathcal{I}_{3}$ &
comb.\\\hline
\multicolumn{1}{|r|}{$5$} & $10.9$ & \multicolumn{1}{|r|}{%
\begin{tabular}
[c]{l}%
$(10.2,11.7)$\\
\multicolumn{1}{r}{$0.431,0.92$}%
\end{tabular}
} & \multicolumn{1}{|r|}{%
\begin{tabular}
[c]{l}%
$(9.9,11.9)$\\
\multicolumn{1}{r}{$0.164,0.95$}%
\end{tabular}
} & \multicolumn{1}{|r|}{%
\begin{tabular}
[c]{l}%
$(10.1,11.7)$\\
\multicolumn{1}{r}{$0.406,0.93$}%
\end{tabular}
} & \multicolumn{1}{|r|}{%
\begin{tabular}
[c]{l}%
$(10.0,11.7)$\\
\multicolumn{1}{r}{$0.93$}%
\end{tabular}
}\\\hline
\multicolumn{1}{|r|}{$10$} & $9.9$ & \multicolumn{1}{|r|}{%
\begin{tabular}
[c]{l}%
$(9.2,10.6)$\\
\multicolumn{1}{r}{$0.176,0.98$}%
\end{tabular}
} & \multicolumn{1}{|r|}{%
\begin{tabular}
[c]{l}%
$(9.3,10.4)$\\
\multicolumn{1}{r}{$0.507,0.93$}%
\end{tabular}
} & \multicolumn{1}{|r|}{%
\begin{tabular}
[c]{l}%
$(9.2,10.5)$\\
\multicolumn{1}{r}{$0.317,0.96$}%
\end{tabular}
} & \multicolumn{1}{|r|}{%
\begin{tabular}
[c]{l}%
$(9.3,10.5)$\\
\multicolumn{1}{r}{$0.95$}%
\end{tabular}
}\\\hline
\multicolumn{1}{|r|}{$25$} & $10.0$ & \multicolumn{1}{|r|}{%
\begin{tabular}
[c]{l}%
$(9.5,10.4)$\\
\multicolumn{1}{r}{$0.192,0.98$}%
\end{tabular}
} & \multicolumn{1}{|r|}{%
\begin{tabular}
[c]{l}%
$(9.6,10.4)$\\
\multicolumn{1}{r}{$0.478,0.96$}%
\end{tabular}
} & \multicolumn{1}{|r|}{%
\begin{tabular}
[c]{l}%
$(9.5,10.4)$\\
\multicolumn{1}{r}{$0.330,0.97$}%
\end{tabular}
} & \multicolumn{1}{|r|}{%
\begin{tabular}
[c]{l}%
$(9.5,10.4)$\\
\multicolumn{1}{r}{$0.97$}%
\end{tabular}
}\\\hline
\multicolumn{1}{|r|}{$100$} & $10.1$ & \multicolumn{1}{|r|}{%
\begin{tabular}
[c]{l}%
$(9.9,10.4)$\\
\multicolumn{1}{r}{$0.229,0.99$}%
\end{tabular}
} & \multicolumn{1}{|r|}{%
\begin{tabular}
[c]{l}%
$(9.9,10.4)$\\
\multicolumn{1}{r}{$0.418,0.99$}%
\end{tabular}
} & \multicolumn{1}{|r|}{%
\begin{tabular}
[c]{l}%
$(9.9,10.4)$\\
\multicolumn{1}{r}{$0.354,0.99$}%
\end{tabular}
} & \multicolumn{1}{|r|}{%
\begin{tabular}
[c]{l}%
$(9.9,10.4)$\\
\multicolumn{1}{r}{$0.99$}%
\end{tabular}
}\\\hline
\end{tabular}
\caption{Relative belief estimates, plausible intervals ( posterior weights, and contents underneath)  for $\mu$  in Example 2.}\label{tab1}%
\end{table}%

Consider now prediction which produces the interesting consequence that
Context II now obtains even when all the models are same.\smallskip

\noindent\textbf{Example 3. }\textit{Prediction.}

Consider Context I but suppose interest is in predicting a future value
$y\in\mathcal{Y},$ whose distribution is conditionally independent of the
observed data $x$ given $\theta$ and has model $\{g_{\lambda}:$ $\lambda
\in\Lambda\}$ where $\Lambda:\Theta\rightarrow\Lambda$ with $\lambda
_{true}=\Lambda(\theta_{true}).$ The first step in solving this problem is to
determine the relevant inference bases and this is done by integrating out the
nuisance parameter which in this case is $\theta.$ So the $i$-th inference
base is given by $\mathcal{I}_{i}=(x,\{m_{i}(\cdot\,|\,y):y\in\mathcal{Y}%
\},m_{i,Y})$ where $m_{i,Y}$ is the density of the $i$-th prior for $y,$
namely, $m_{i,Y}(y)=\int_{\Theta}g_{\Lambda(\theta)}(y)\pi_{i}(\theta
)\,d\theta,$ and $m_{i}(x\,|\,y)=\int_{\Theta}f_{\theta}(x)g_{\Lambda(\theta
)}(y)\pi_{i}(\theta)\,d\theta/m_{i,Y}(y)$ is the conditional density of $x$
given $y.$ Note that unconditionally $x$ and $y$ are not independent and now
the collection of possible distributions for $x$ is indexed by $y.$ The $i$-th
posterior density of $y$ is then $m_{i,Y}(y\,|\,x)=m_{i}(x\,|\,y)m_{i,Y}%
(y)/m_{i}(x).$

The models $\{m_{i}(\cdot\,|\,y):y\in\mathcal{Y}\}$ are now not all the same
so this is Context II. It is assumed, as is typically the case, that the mss
for these models is constant in $i$ so the weights are comparable. Applying
(\ref{linrule2}), with the single data set $x,$ leads to
\[
RB_{1,\alpha,Y}^{\ast}(y\,|\,x)=\sum_{i=1}^{k}\frac{\alpha_{i}m_{i}%
(x)}{m_{1,\alpha}(x)}RB_{i,Y}(y\,|\,x)
\]
with $RB_{i,Y}(y\,|\,x)=m_{i,Y}(y\,|\,x)/m_{i,Y}(y)=m_{i}(x\,|\,y)/m_{i}(x)$
and (\ref{Jeffreypost}) leads to posterior $m_{1,\alpha,Y}^{\ast
}(y\,|\,x)=\sum_{i=1}^{k}(\alpha_{i}m_{i}(x)/m_{1,\alpha}(x))m_{i,Y}%
(y\,|\,x).$ Note that in this case\ the posterior of $y$ given $x$ is
well-defined via Bayesian conditioning and equals $m_{1,\alpha,Y}^{\ast
}(y\,|\,x)\ $so there is no need to invoke Jeffrey's conditionalization for
the posterior. It is notable, however, that if the relative belief ratio for
$y$ is computed using this posterior and the prior $m_{1,\alpha,Y}%
(y)=\sum_{i=1}^{k}\alpha_{i}m_{i,Y}(y),$ then this equals%
\begin{equation}
\sum_{i=1}^{k}\frac{\alpha_{i}m_{i}(x)m_{i,Y}(y)}{m_{1,\alpha}(x)m_{1,\alpha
,Y}(y)}RB_{i,Y}(y\,|\,x) \label{rbfullpred}%
\end{equation}
which does not equal $RB_{1,\alpha,Y}^{\ast}(y\,|\,x).$ Given that the weights
in (\ref{rbfullpred}) depend on the object of interest $y,$ this does not
correspond to linear pooling of the evidence and this is because the model is
not constant. There is no reason to suppose that (\ref{rbfullpred}) will
retain the good properties of linear pooling and experience with it suggests
that it is not the correct way to combine. As such, the recommended approach
is via (\ref{linrule2}) based on Jeffrey's conditionization and which retains
the good properties of linear pooling.

Suppose now the context is as discussed in Example 2 but the goal is to make a
prediction concerning a future independent value $y\sim N(\mu,\sigma_{0}%
^{2}).$ So the $i$-th prior $m_{i,Y}$ is given by $y\sim N(\mu_{i},\sigma
_{0}^{2}+\tau_{i}^{2})$ and the $i$-th posterior $m_{i,Y}(\cdot\,|\,\bar{x})$
is $y\,|\,\bar{x}\sim N((n/\sigma_{0}^{2}+1/\tau_{i}^{2})^{-1}(n\bar{x}%
/\sigma_{0}^{2}+\mu_{i}/\tau_{i}^{2}),(n/\sigma_{0}^{2}+1/\tau_{i}^{2}%
)^{-1}+\tau_{i}^{2}).$ Table 2 gives the results for predicting $y$ using the
data in Example 2.\ The final row indicates what happens as $n\rightarrow
\infty$ and note that the weights converge as well with the $i$-th limiting
weight proportional to $(\sigma_{0}^{2}+\tau_{i}^{2})^{-1/2}\exp(-(\mu-\mu
_{i})^{2}/2(\sigma_{0}^{2}+\tau_{i}^{2}))$ which depends on the relative
accuracy of the $i$-th prior with respect to the true mean $\mu.$ When all the
prior variances are the same the prior which has its mean closest to the true
value will give the heaviest weight. Also, as $\tau_{i}^{2}\rightarrow\infty$
the $i$-th weight goes to 0. Note that the limiting plausible intervals are
dependent on the prior and the interval does not shrink to a point because $y$
is random. The limiting posterior content of these intervals is the
probability content given by the true distribution of $y.$

For the limiting plausible intervals for $y$ to still be dependent on the
prior is different than the situation when making inference about a parameter
as, in that case, the plausible intervals shrink to the true value as the
amount of data increases. The difference is that there is not a "true" value
for $y.$ The limiting plausible interval does not allow for all possible
values for $y$ and the effect of the prior is to disallow some possible values
because belief in such a value is less than that specified by the prior of
$y.$ As can be seen from Table 2 this effect is not great unless the prior, as
with $\pi_{1}$ here, puts little mass near the true value. However, such an
occurrence also reduces the limiting weight for such a component.
$\blacksquare$\smallskip%

\begin{table}[tbp] \centering
\begin{tabular}
[c]{|c|c|c|c|c|c|}\hline
$n$ & pred. & $\mathcal{I}_{1}$ & $\mathcal{I}_{2}$ & $\mathcal{I}_{3}$ &
comb.\\\hline
\multicolumn{1}{|r|}{$5$} & $10.9$ & \multicolumn{1}{|r|}{%
\begin{tabular}
[c]{l}%
$(8.7,12.1)$\\
\multicolumn{1}{r}{$0.82$}%
\end{tabular}
} & \multicolumn{1}{|r|}{%
\begin{tabular}
[c]{l}%
$(9.7,16.0)$\\
\multicolumn{1}{r}{$0.81$}%
\end{tabular}
} & \multicolumn{1}{|r|}{%
\begin{tabular}
[c]{l}%
$(9.4,12.4)$\\
\multicolumn{1}{r}{$0.83$}%
\end{tabular}
} & \multicolumn{1}{|r|}{%
\begin{tabular}
[c]{l}%
$(9.2,13.7)$\\
\multicolumn{1}{r}{$0.94$}%
\end{tabular}
}\\\hline
\multicolumn{1}{|r|}{$10$} & $9.9$ & \multicolumn{1}{|r|}{%
\begin{tabular}
[c]{l}%
$(6.5,11.1)$\\
\multicolumn{1}{r}{$0.87$}%
\end{tabular}
} & \multicolumn{1}{|r|}{%
\begin{tabular}
[c]{l}%
$(9.03,12.5)$\\
\multicolumn{1}{r}{$0.79$}%
\end{tabular}
} & \multicolumn{1}{|r|}{%
\begin{tabular}
[c]{l}%
$(8.0,11.2)$\\
\multicolumn{1}{r}{$0.86$}%
\end{tabular}
} & \multicolumn{1}{|r|}{%
\begin{tabular}
[c]{l}%
$(7.7,11.8)$\\
\multicolumn{1}{r}{$0.94$}%
\end{tabular}
}\\\hline
\multicolumn{1}{|r|}{$25$} & $10.0$ & \multicolumn{1}{|r|}{%
\begin{tabular}
[c]{l}%
$(6.7,11.2)$\\
\multicolumn{1}{r}{$0.87$}%
\end{tabular}
} & \multicolumn{1}{|r|}{%
\begin{tabular}
[c]{l}%
$(9.1,12.7)$\\
\multicolumn{1}{r}{$0.79$}%
\end{tabular}
} & \multicolumn{1}{|r|}{%
\begin{tabular}
[c]{l}%
$(8.2,11.3)$\\
\multicolumn{1}{r}{$0.86$}%
\end{tabular}
} & \multicolumn{1}{|r|}{%
\begin{tabular}
[c]{l}%
$(7.9,11.9)$\\
\multicolumn{1}{r}{$0.95$}%
\end{tabular}
}\\\hline
\multicolumn{1}{|r|}{$100$} & $10.1$ & \multicolumn{1}{|r|}{%
\begin{tabular}
[c]{l}%
$(7.1,11.3)$\\
\multicolumn{1}{r}{$0.87$}%
\end{tabular}
} & \multicolumn{1}{|r|}{%
\begin{tabular}
[c]{l}%
$(9.3,13.2)$\\
\multicolumn{1}{r}{$0.80$}%
\end{tabular}
} & \multicolumn{1}{|r|}{%
\begin{tabular}
[c]{l}%
$(8.4,11.4)$\\
\multicolumn{1}{r}{$0.86$}%
\end{tabular}
} & \multicolumn{1}{|r|}{%
\begin{tabular}
[c]{l}%
$(8.1,12.2)$\\
\multicolumn{1}{r}{$0.96$}%
\end{tabular}
}\\\hline
$\infty$ & $10.0$ &
\begin{tabular}
[c]{l}%
$(6.9,11.2)$\\
\multicolumn{1}{r}{$0.88$}%
\end{tabular}
&
\begin{tabular}
[c]{l}%
$(9.2,12.8)$\\
\multicolumn{1}{r}{$0.79$}%
\end{tabular}
&
\begin{tabular}
[c]{l}%
$(8.2,11.3)$\\
\multicolumn{1}{r}{$0.87$}%
\end{tabular}
&
\begin{tabular}
[c]{l}%
$(8.0,12.0)$\\
\multicolumn{1}{r}{$0.96$}%
\end{tabular}
\\\hline
\end{tabular}
\caption{Relative belief estimates, plausible intervals (posterior contents underneath) for $y$  in Example 3 with the posterior weights  as in Table 1.}\label{tab1 copy(1)}%
\end{table}%

Consider now an example where the weights require adjustment.\smallskip

\noindent\textbf{Example 4. }\textit{Location-normal models with different
variances.}

Suppose a situation similar to Example 2 but now with three distinct models so
this is Context II. Here the $i$-th statistician assumes that the true
distribution is $N(\mu,\sigma_{i0}^{2})$ where the $\sigma_{i0}^{2}$ are known
but $\mu\in%
\mathbb{R}
^{1}$ is unknown and interest is in $\psi=\Psi(\mu)=\mu.$ Again three
$N(\mu_{i0},\tau_{i0}^{2})$ priors are assumed. So the statisticians disagree
about the "known" variance of the sampling distribution and an ancillary needs
to play a role to make the weights comparable.

In this case $A(x)=x-\bar{x}1$ is ancillary for each model and is
independently distributed from the common mss $L(x)=\bar{x}\sim N(\mu
,\sigma_{i0}^{2}/n)$ and $x\leftrightarrow(L(x),A(x)).$ Therefore, with equal
weights for the priors, and taking the ancillaries into account, the $i$-th
weight satisfies%
\[
m_{i}(L(x)\,|\,A(x))\propto(\sigma_{i0}^{2}/n+\tau_{i}^{2})^{-1/2}%
\varphi\left(  (\sigma_{i0}^{2}/n+\tau_{i}^{2})^{-1/2}\bar{x}\right)  .
\]
From this it is seen that the assumed variances and the prior both play a role
in determining how much weight a given analysis should have. Note that as
$\sigma_{i0}^{2}\rightarrow\infty$ or $\tau_{i}^{2}\rightarrow\infty$, and all
other parameters are fixed, then the weight of the $i$-th analysis goes to 0
as it should as, in the limit, no information is being provided about the true
value of $\mu.$ Proposition 8 tell us that when $n\rightarrow\infty$ and the
$i$-th variance is correct and the others are not, then the $i$-th inference
base will dominate. $\blacksquare$\smallskip

Consider now an example where the models are truly different.\smallskip

\noindent\textbf{Example 5. }\textit{Location with quite different models.}

Consider now the context of Example 2 and suppose that one of the models, say
the one in $\mathcal{I}_{1},$ is a $t_{1}$ (Cauchy) location model, while the
other models and all the priors are as previously specified. For all three
inference bases $A(x)=x-\bar{x}1$ is ancillary. To insure that $\sigma_{0}$
has the same interpretation across all inference bases, the $t_{1}$ density is
rescaled by $\eta_{0}$ so that the interval $(-\sigma_{0},\sigma_{0})$
contains $0.6827$ of the probability for all 3 distributions. This implies
$\eta_{0}=\sigma_{0}/\tan(0.1827\pi)$ and, with $g(z)=1/\eta_{0}\pi
(1+z^{2}/\eta_{0}^{2}),$ the first model is $\{f_{\mu}(x):\mu\in%
\mathbb{R}
^{1}\}$ where $f_{\mu}(x)=\prod\nolimits_{i=1}^{n}g(x_{i}-\mu).$ To obtain the
corresponding weight the following expression needs to be evaluated numerically,%

\[
m_{1}(x\,|\,A(x))=\frac{\int_{-\infty}^{\infty}f_{1,\mu}((\bar{x}%
-\mu)1+A(x))\pi_{1}(\mu)\,d\mu}{\int_{-\infty}^{\infty}\int_{-\infty}^{\infty
}f_{1,\mu}((\bar{x}-\mu)1+A(x))\pi_{1}(\mu)\,d\mu\,d\bar{x}}.
\]
When applied to the data of Example 2 very similar results are obtained. Table
3 contains the weights for the inference bases for this situation.
$\blacksquare$\smallskip%

\begin{table}[tbp] \centering
\begin{tabular}
[c]{|c|c|c|c|}\hline
$n$ & $\mathcal{I}_{1}$ & $\mathcal{I}_{2}$ & $\mathcal{I}_{3}$\\\hline
\multicolumn{1}{|r|}{$5$} & \multicolumn{1}{|r|}{$0.4248$} &
\multicolumn{1}{|r|}{$0.1652$} & \multicolumn{1}{|r|}{$0.4100$}\\\hline
\multicolumn{1}{|r|}{$10$} & \multicolumn{1}{|r|}{$0.2036$} &
\multicolumn{1}{|r|}{$0.4900$} & \multicolumn{1}{|r|}{$0.3064$}\\\hline
\multicolumn{1}{|r|}{$25$} & \multicolumn{1}{|r|}{$0.1716$} &
\multicolumn{1}{|r|}{$0.4898$} & \multicolumn{1}{|r|}{$0.3386$}\\\hline
\multicolumn{1}{|r|}{$100$} & \multicolumn{1}{|r|}{$0.2235$} &
\multicolumn{1}{|r|}{$0.4202$} & \multicolumn{1}{|r|}{$0.3563$}\\\hline
\end{tabular}
\caption{Weights for the inference bases  in Example 5}\label{tab1 copy(2)}%
\end{table}%

Consider now a common context for applications.\textbf{ }\smallskip

\noindent\textbf{Example 6. }\textit{Linear regression.}

Suppose that the data is $(x_{i},y_{i})$ for $i=1,\ldots,n$ and there are two
analysts where both propose a simple regression model $y\,=X\beta+\sigma z$
where $X=(1_{n}/\sqrt{n},x)$ with $1_{n}\perp x$ and $||x||=1,\beta=(\beta
_{1},\beta_{2})^{\prime}\in%
\mathbb{R}
^{2}$ and $\sigma>0$ unknown and $z$ is a sample from $N(0,1)$ for analyst 1
and is a sample from a $t_{\lambda}/\sqrt{(\lambda-2)/\lambda}$
distribution\ for analyst 2 for some value $\lambda>2.$ In both models
$\sigma^{2}$ is the variance of a $y_{i}.$ Letting $b=(X^{\prime}%
X)^{-1}X^{\prime}y$ be the least squares estimate of $\beta$ and
$s^{2}=||y-Xb||,$ then $y\leftrightarrow(L(y),A(y))$ where $L(y)=(b,s^{2})$
and $A(y)=(y-Xb)/s$ is ancillary for both models. Further suppose that the
quantity of inferential interest is the slope parameter $\psi=\Psi(\beta
_{1},\beta_{2},\sigma^{2})=\beta_{2}.$ Denoting the relevant density of a
$z_{i}$ by $f,$ the joint density of $(b,s)$ given $A(y)=a$ is proportional
to
\[
s^{n-3}\sigma^{-n}\prod\nolimits_{i=1}^{n}f\left(  \frac{b_{1}-\beta_{1}%
}{\sigma}+\frac{b_{2}-\beta_{2}}{\sigma}x_{i}+\frac{s}{\sigma}a_{i}\right)  .
\]
The posterior density of $\beta_{2}$ can be worked out in closed-form when $f$
is the $N(0,1)$ density but generally it will require numerical integration to
determine the posterior density and the posterior weights for the combination.

For the prior suppose both analysts agree on $\beta\,|\,\sigma^{2}\sim
N_{2}(0,\tau_{0}^{2}\sigma^{2}I)$ and $1/\sigma^{2}\sim$ gamma$_{\text{rate}%
}(\alpha_{1},\alpha_{2}).$ Note that the zero mean for $\beta$ may entail
subtracting a known, fixed constant vector from\ $y$ so this, and the
assumption that $1_{n}\perp x,$ may entail some preprocessing of the data. The
prior distribution of the quantity of interest is then $\beta_{2}\sim\tau
_{0}\sqrt{\alpha_{2}/\alpha_{1}}t_{2\alpha_{1}}$ where $t_{2\alpha_{1}}$
denotes the $t$ distribution on $2\alpha_{1}$ degrees of freedom.

To obtain the hyperparameters of the prior requires elicitation and this can
be carried out using the following method as described in Evans and Tomal
(2018). Suppose that it is known with virtual certainty, based on our
knowledge of the measurements being taken, that $\beta_{1}+\beta_{2}x$ will
lie in the interval $(-m_{0},m_{0})$ for some $m_{0}>0$ for all $x\in R$ a
compact set centered at 0 and contained in $[-1,1]^{k}$ on account of the
standardization. The phrase `virtual certainty' is interpreted here as a
probability greater than or equal to $\gamma$ where $\gamma$ is some large
probability like $0.99.$ Therefore, the prior on $\beta$ must satisfy
$2\Phi(m_{0}/\sigma\tau_{0}\{1+x^{2}\}^{1/2})-1\geq\gamma$ for all $x\in R$
which implies%
\begin{equation}
\sigma\leq m_{0}/\zeta_{0}\tau_{0}z_{(1+\gamma)/2} \label{ineq1}%
\end{equation}
where $\zeta_{0}^{2}=1+\max_{x\in R}x^{2}\leq2$ with equality when $R=[-1,1].$
An interval that will contain a response value $y$ with virtual certainty,
given predictor value $x,$ is $\beta_{1}+\beta_{2}x\pm\sigma z_{(1+\gamma
)/2}.$ Suppose that we have lower and upper bounds $s_{1}$ and $s_{2}$ on the
half-length of this interval so that $s_{1}\leq\sigma z_{(1+\gamma)/2}\leq
s_{2}$ or, equivalently,%
\begin{equation}
s_{1}/z_{(1+\gamma)/2}\leq\sigma\leq s_{2}/z_{(1+\gamma)/2} \label{ineq2}%
\end{equation}
holds with virtual certainty. Combining (\ref{ineq2}) with (\ref{ineq1})
implies $\tau_{0}=m_{0}/s_{2}\zeta_{0}.$ To obtain the relevant values of
$\alpha_{1}$ and $\alpha_{2}$ let $G\left(  \alpha_{1},\alpha_{2}%
,\cdot\right)  $ denote the cdf of the gamma$_{\text{rate}}\left(  \alpha
_{1},\alpha_{2}\right)  $ distribution and note that $G\left(  \alpha
_{1},\alpha_{2},z\right)  =G\left(  \alpha_{1},1,\alpha_{2}z\right)  .$
Therefore, the interval for $1/\sigma^{2}$ implied by (\ref{ineq2}) contains
$1/\sigma^{2}$ with virtual certainty, when $\alpha_{1},\alpha_{2}$ satisfy
$G^{-1}(\alpha_{1},\alpha_{2},(1+\gamma)/2)=s_{1}^{-2}z_{(1+\gamma)/2}%
^{2},G^{-1}(\alpha_{1},\alpha_{2},(1-\gamma)/2)=s_{2}^{-2}z_{(1-\gamma)/2}%
^{2},$ or equivalently
\begin{align}
G(\alpha_{1},1,\alpha_{2}s_{1}^{-2}z_{(1+\gamma)/2}^{2})  &  =(1+\gamma
)/2,\label{eq3}\\
G(\alpha_{1},1,\alpha_{2}s_{2}^{-2}z_{(1-\gamma)/2}^{2})  &  =(1-\gamma)/2.
\label{eq4}%
\end{align}
It is a simple matter to solve these equations for $\left(  \alpha_{1}%
,\alpha_{2}\right)  .$ For this choose an initial value for $\alpha_{1}$ and,
using (\ref{eq3}), find $z$ such that $G(\alpha_{1},1,z)=(1+\gamma)/2,$ which
implies $\alpha_{2}=zs_{1}^{2}/z_{(1+\gamma)/2}^{2}.$ If the left-side of
(\ref{eq4}) is less (greater) than $(1-\gamma)/2,$ then decrease
(increase)\ the value of $\alpha_{1}$ and repeat step 1. Continue iterating
this process until satisfactory convergence is attained.

Consider now a numerical example drawn from Zellner (1996) where the response
variable is income in U.S. dollars per capita (deflated), and the predictor
variable is investment in dollars per capita (deflated) for the United States
for the years 1922--1941. The data are provided in Table \ref{zellner}. The
data vector $y$ was replaced by $y-X(340,3)^{t}$ as this centered the
observations about 0. Taking $\gamma=0.99,\zeta_{0}=\sqrt{2},m_{0}%
=30,s_{1}=10,s_{2}=40$ leads to the values $\tau_{0}=0.54,\alpha
_{1}=4.05,\alpha_{2}=140.39.$ The following prior is then used for both
models,
\[
(\beta_{1},\beta_{2})\,|\,,\sigma^{2}\sim N_{2}(0,(0.54)^{2}\sigma
^{2}I),\text{ }1/\sigma^{2}\sim\text{gamma}(4.05,140.39).
\]
Table \ref{weightsZel} presents the weights that result when different
$t_{\lambda}$ error distributions are considered to be combined with the
results from a $N(0,1)$ error assumption. Presumably this arises when one
analyst is concerned that tails longer than the normal are appropriate. As can
be seen the normal error assumption dominates except for $\lambda=100$ when
the inferences don't differ by much in any case. This is not surprising as
various residual plots don't indicate any issue with the normality assumption
for these data. These weights were computed using importance sampling and were
found to be robust to the prior by repeating the computations after making
small changes to the hyperparameters.

The approach taken in this example is easily generalized to more general
linear regression models including situations where the priors change.
$\blacksquare$%

\begin{table}[tbp] \centering
$%
\begin{tabular}
[c]{|llllll|}\hline
Year & Income & \multicolumn{1}{l|}{Investment} & Year & Income &
Investment\\\hline
\multicolumn{1}{|r}{$1922$} & \multicolumn{1}{r}{$433$} &
\multicolumn{1}{r|}{$39$} & \multicolumn{1}{r}{$1932$} &
\multicolumn{1}{r}{$372$} & \multicolumn{1}{r|}{$22$}\\
\multicolumn{1}{|r}{$1923$} & \multicolumn{1}{r}{$483$} &
\multicolumn{1}{r|}{$60$} & \multicolumn{1}{r}{$1933$} &
\multicolumn{1}{r}{$381$} & \multicolumn{1}{r|}{$17$}\\
\multicolumn{1}{|r}{$1924$} & \multicolumn{1}{r}{$479$} &
\multicolumn{1}{r|}{$42$} & \multicolumn{1}{r}{$1934$} &
\multicolumn{1}{r}{$419$} & \multicolumn{1}{r|}{$27$}\\
\multicolumn{1}{|r}{$1925$} & \multicolumn{1}{r}{$486$} &
\multicolumn{1}{r|}{$52$} & \multicolumn{1}{r}{$1935$} &
\multicolumn{1}{r}{$449$} & \multicolumn{1}{r|}{$33$}\\
\multicolumn{1}{|r}{$1926$} & \multicolumn{1}{r}{$494$} &
\multicolumn{1}{r|}{$47$} & \multicolumn{1}{r}{$1936$} &
\multicolumn{1}{r}{$511$} & \multicolumn{1}{r|}{$48$}\\
\multicolumn{1}{|r}{$1927$} & \multicolumn{1}{r}{$498$} &
\multicolumn{1}{r|}{$51$} & \multicolumn{1}{r}{$1937$} &
\multicolumn{1}{r}{$520$} & \multicolumn{1}{r|}{$51$}\\
\multicolumn{1}{|r}{$1928$} & \multicolumn{1}{r}{$511$} &
\multicolumn{1}{r|}{$45$} & \multicolumn{1}{r}{$1938$} &
\multicolumn{1}{r}{$477$} & \multicolumn{1}{r|}{$33$}\\
\multicolumn{1}{|r}{$1929$} & \multicolumn{1}{r}{$534$} &
\multicolumn{1}{r|}{$60$} & \multicolumn{1}{r}{$1939$} &
\multicolumn{1}{r}{$517$} & \multicolumn{1}{r|}{$46$}\\
\multicolumn{1}{|r}{$1930$} & \multicolumn{1}{r}{$478$} &
\multicolumn{1}{r|}{$39$} & \multicolumn{1}{r}{$1940$} &
\multicolumn{1}{r}{$548$} & \multicolumn{1}{r|}{$54$}\\
\multicolumn{1}{|r}{$1931$} & \multicolumn{1}{r}{$440$} &
\multicolumn{1}{r|}{$41$} & \multicolumn{1}{r}{$1941$} &
\multicolumn{1}{r}{$629$} & \multicolumn{1}{r|}{$100$}\\\hline
\end{tabular}
\ $%
\caption{Haavelmo's data on income and investment from Zellner (1996).}\label{zellner}%
\end{table}%
%

\begin{table}[tbp] \centering
\begin{tabular}
[c]{|l|c|c|c|c|c|c|}\hline
$\lambda$ & $100$ & $50$ & $20$ & $10$ & $5$ & $3$\\\hline
$N(0,1)$ & $0.556$ & $0.612$ & $0.766$ & $0.928$ & $0.998$ & $1.000$\\
$t_{\lambda}$ & $0.444$ & $0.388$ & $0.234$ & $0.072$ & $0.002$ &
$0.000$\\\hline
\end{tabular}
\caption{Weights for normal and $t_\lambda$ errors in Example 6.}\label{weightsZel}%
\end{table}%

\section{Conclusions}

The problem of how to combine evidence has been considered for a Bayesian
context where each analyst proposes a model and prior for the same data.
Linear opinion pooling is seen as the natural way to make such a combination
at least when the inference bases only differ in the priors on the parameter
of interest. This has been shown to have appropriate properties such as
preserving a consensus with respect to the evidence and, when combining
evidence is considered as opposed to just combining priors, behaves
appropriately when considering independent events. In certain contexts the
idea can be extended in a logical way based on the idea underlying Jeffrey
conditionalization. There are restrictions as in the end the posterior weights
have to be seen to be comparable and focused on that aspect of the data which
is relevant for inference about the unknowns. Asymptotically the approach has
been shown to behave correctly.

Certainly this does not cover all contexts where one might want to combine
evidence as when there are different data sets and different models. If the
models are all for the same response variable, then one possibility is to
simply combine data sets and proceed as we have done here. More generally it
may be that the only aspect in common among the models is the characteristic
of interest $\Psi$ and then it is not clear how we should combine and this
warrants further investigation.

\section{Appendix}

\subparagraph{\noindent Proof of Proposition 7}

Parts (i) and (ii) are established in Evans (2015), Section 4.7 for a general
prior and so can be applied with the prior $\pi_{1,\alpha}.$ For part (iii) we
have
\begin{align*}
&  \frac{\alpha_{i}m_{i}(x)}{m_{1,\alpha}(x)}=\frac{\alpha_{i}\sum_{\theta
}f_{\theta}(\breve{x}_{n})\pi_{i}(\theta)}{\sum_{\theta}f_{\theta}(\breve
{x}_{n})\pi_{j}(\theta)}\pi_{1,\alpha}(\theta)\\
&  =\frac{\alpha_{i}\sum_{\theta}\exp\left\{  -n\left(  \frac{1}{n}\log
\frac{f_{\theta_{true}}(\breve{x}_{n})}{f_{\theta}(\breve{x}_{n})}\right)
\right\}  \pi_{i}(\theta)}{\sum_{\theta}\exp\left\{  -n\left(  \frac{1}{n}%
\log\frac{f_{\theta_{true}}(\breve{x}_{n})}{f_{\theta}(\breve{x}_{n})}\right)
\right\}  \pi_{1,\alpha}(\theta)}\\
&  =\frac{\alpha_{i}\left[  \pi_{i}(\theta_{true})+\sum_{\theta\neq
\theta_{true}}\exp\left\{  -n\left(  \frac{1}{n}\log\frac{f_{\theta_{true}%
}(\breve{x}_{n})}{f_{\theta}(\breve{x}_{n})}\right)  \right\}  \pi_{i}%
(\theta)\right]  }{\pi_{1,\alpha}(\theta_{true})+\sum_{\theta\neq\theta
_{true}}\exp\left\{  -n\left(  \frac{1}{n}\log\frac{f_{\theta_{true}}%
(\breve{x}_{n})}{f_{\theta}(\breve{x}_{n})}\right)  \right\}  \pi_{1,\alpha
}(\theta)}%
\end{align*}
and by the SLLN $\frac{1}{n}\log f_{\theta_{true}}(\breve{x}_{n})/f_{\theta
}(\breve{x}_{n})\rightarrow KL(f_{\theta_{true}},f_{\theta}),$ where $KL$ is
the Kullback-Leibler divergence. Since $KL(f_{\theta_{true}},f_{\theta})\geq0$
and 0 iff $\theta=\theta_{true},$ this completes the result. Part (iv) is
established similarly. $\blacksquare$

\subparagraph{Proof of Proposition 8}

Suppose initially that only one of the proposed models contains the true
distribution and wlog it is given by $i=1.$ Following the proof of Proposition
7 (iii) then
\begin{align*}
&  \frac{\alpha_{i}m_{i}(\breve{x}_{n})}{m_{1,\alpha}(\breve{x}_{n})}%
=\frac{\alpha_{i}\sum_{\theta_{i}\in\Theta_{i}}f_{\theta_{i}}(\breve{x}%
_{n})\pi_{i}(\theta_{i})}{\sum_{j}\alpha_{j}\sum_{\theta_{j}\in\Theta_{j}%
}f_{\theta_{j}}(\breve{x}_{n})\pi_{j}(\theta_{j})}\\
&  =\frac{\alpha_{i}\sum_{\theta_{i}\in\Theta_{i}}\exp\left\{  -n\left(
\frac{1}{n}\log\frac{f_{1\theta_{1true}}(\breve{x}_{n})}{f_{\theta_{i}}%
(\breve{x}_{n})}\right)  \right\}  \pi_{i}(\theta_{i})}{\sum_{j}\alpha_{j}%
\sum_{\theta_{j}\in\Theta_{j}}\exp\left\{  -n\left(  \frac{1}{n}\log
\frac{f_{1\theta_{1true}}(\breve{x}_{n})}{f_{\theta_{j}}(\breve{x}_{n}%
)}\right)  \right\}  \pi_{j}(\theta_{j})}\rightarrow\left\{
\begin{array}
[c]{cc}%
1 & i=1,\\
0 & i\neq1.
\end{array}
\right.
\end{align*}
If two of the models contain the true distribution, say given by $i=1,2,$
then
\[
\frac{\alpha_{i}m_{i}(\breve{x}_{n})}{m_{1,\alpha}(\breve{x}_{n})}%
\rightarrow\left\{
\begin{array}
[c]{cc}%
\frac{\alpha_{i}\pi_{i}(\theta_{itrue})}{\alpha_{1}\pi_{1}(\theta
_{1true})+\alpha_{2}\pi_{2}(\theta_{2true})} & i=1,2,\\
0 & i\neq1,2.
\end{array}
\right.
\]
This line of reasoning proves (i).

Now note that%
\begin{align*}
&  RB_{i,\Psi}(\psi\,|\,\breve{x}_{n})=\frac{m_{i}(\breve{x}_{n}\,|\,\psi
)}{m_{i}(\breve{x}_{n})}=\frac{\sum_{\theta_{i}\in\Psi^{-1}\{\psi\}}%
f_{i\theta_{i}}(\breve{x}_{n})\pi_{i}(\theta_{i}\,|\,\psi)}{\sum_{\theta
_{i}\in\Theta_{i}}f_{i\theta_{i}}(\breve{x}_{n})\pi_{i}(\theta_{i}\,|\,\psi
)}\\
&  =\frac{1}{\pi_{i,\Psi}(\psi)}\frac{\sum_{\theta_{i}\in\Psi^{-1}\{\psi
\}}\exp\left\{  -n\left(  \frac{1}{n}\log\frac{f_{1\theta_{1true}}(\breve
{x}_{n})}{f_{\theta_{i}}(\breve{x}_{n})}\right)  \right\}  \pi_{i}(\theta
_{i})}{\sum_{\theta_{i}\in\Theta_{i}}\exp\left\{  -n\left(  \frac{1}{n}%
\log\frac{f_{1\theta_{1true}}(\breve{x}_{n})}{f_{\theta_{i}}(\breve{x}_{n}%
)}\right)  \right\}  \pi_{i}(\theta_{i})}%
\end{align*}
which implies $0<RB_{i,\Psi}(\psi\,|\,\breve{x}_{n})\leq1/\pi_{i,\Psi}(\psi).$
If $i\notin J,$ then the $i$-ith term in $RB_{1,\alpha,\Psi}^{\ast}%
(\psi\,|\,\breve{x}_{n})$ converges to 0 as it has been shown that the $i$-th
weight does. If $i\in J,$ then $RB_{i,\Psi}(\psi\,|\,\breve{x}_{n})\rightarrow
I_{\{\psi_{true}\}}(\psi)/\pi_{i\Psi}(\psi)$ which proves the first part of
(ii). Now note that $\sum_{i\in J}w_{i}/\pi_{i,\Psi}(\psi_{true})\geq
1/\max\{\pi_{i,\Psi}(\psi)\}>1$ proving the second part.

Now $\pi_{i,\Psi}(\psi\,|\,x)$ is bounded so it is only necessary to the limit
when $i\in J$ and in that case $\pi_{i,\Psi}(\psi\,|\,x)\rightarrow
I_{\{\psi_{true}\}}(\psi)$ and the result follows. $\blacksquare$

\subparagraph{Proof of Proposition 9}

As in the proof of Proposition 7, suppose that the true distribution is
contained in only one of the models, say $i=1.$%
\begin{align*}
&  \frac{\alpha_{i}m_{i}(x\,|\,A(x))}{m_{1,\alpha}(\breve{x}_{n}%
\,|\,A(x))}=\frac{\alpha_{i}\sum_{\theta_{i}\in\Theta_{i}}(f_{\theta_{i}%
}(\breve{x}_{n})/f_{i}(\breve{n}(\breve{x}_{n})))\pi_{i}(\theta_{i})}{\sum
_{j}\alpha_{j}\sum_{\theta_{j}\in\Theta_{j}}(f_{\theta_{j}}(\breve{x}%
_{n})/(\breve{n}(\breve{x}_{n})))/\pi_{j}(\theta_{j})}\\
&  =\frac{\alpha_{i}\sum_{\theta_{i}\in\Theta_{i}}\exp\left\{  -n\left(
\frac{1}{n}\log\frac{f_{1\theta_{1true}}(\breve{x}_{n})/f_{1}(\breve{n}%
(\breve{x}_{n}))}{f_{\theta_{i}}(\breve{x}_{n})/f_{i}(\breve{n}(\breve{x}%
_{n}))}\right)  \right\}  \pi_{i}(\theta_{i})}{\sum_{j}\alpha_{j}\sum
_{\theta_{j}\in\Theta_{j}}\exp\left\{  -n\left(  \frac{1}{n}\log
\frac{f_{1\theta_{1true}}(\breve{x}_{n})/f_{1}(\breve{n}(\breve{x}_{n}%
)}{f_{\theta_{j}}(\breve{x}_{n})/f_{j}(\breve{n}(\breve{x}_{n}))}\right)
\right\}  \pi_{j}(\theta_{j})}%
\end{align*}
Now noting that
\begin{align*}
n\log\frac{f_{j}(\breve{n}(\breve{x}_{n}))}{f_{1}(\breve{n}(\breve{x}_{n}))}
&  =\log\left\{  \left(  \frac{p_{j1}}{p_{11}}\right)  ^{n_{1}(\breve{x}%
_{n})/n}\cdots\left(  \frac{p_{jm}}{p_{1m}}\right)  ^{n_{m}(\breve{x}_{n}%
)/n}\right\} \\
&  \rightarrow\log\left\{  \left(  \frac{p_{j1}}{p_{11}}\right)  ^{p_{11}%
}\cdots\left(  \frac{p_{jm}}{p_{1m}}\right)  ^{p_{1m}}\right\}
\end{align*}
the result is obtained. The remainder of the proof is as in Proposition 8.
$\blacksquare$

\section{References}

\noindent Albert, I., Donnet, S., Guihenneuc-Jouyaux, C., Low-Choy S.,
Mengersen, K. and Rousseau, J. (2012) Combining expert opinions in prior
elicitation. Bayesian Analysis, 7, 3, 503-532.\smallskip

\noindent Birnbaum, A. (1962) On the foundations of statistical inference. J.
of the American Statistical Association, 57, 298, 269-306.\smallskip

\noindent Bunn, D. K. (1981) Two methodologies for the linear combination of
forecasts. J. of the Operational Research Society, 32, 213-222.\smallskip

\noindent Burgman, M. A., McBride , M., Ashton, R., Speirs-Bridge, A.,
Flander, L., Wintle, B, Rumpff, L. and Twardy, C. (2011) Expert status and
performance. PLoS ONE, 6, e22998, doi:10.1371/journal.pone.0022998.\smallskip

\noindent Clemen, R. T. and Winkler, R. L. (1999) Combining probability
distributions from experts in risk analysis. Risk Analysis, 19 (2),
187-203.\smallskip\ 

\noindent Cooke, R. M. (1991) Experts in Uncertainty: Opinion and Subjective
Probability in Science. Oxford University Press.\smallskip\ 

\noindent DeGroot, M. H. (1974) Reaching a consensus. J. of the American
Statistical Association, 69, 118-121.\smallskip

\noindent Diaconis, P. and Zabell, S. (1982) Updating subjective probability.
Journal of the American Statistical Association, 77, 380, 822-820\smallskip

\noindent Evans, M. (2015) Measuring Statistical Evidence Using Relative
Belief. Monographs on Statistics and Applied Probability 144, CRC Press,
Taylor \& Francis Group.\smallskip

\noindent Evans, M. and Jang, G-H. (2011). Weak informativity and the
information in one prior relative to another. Statistical Science, 26, 3,
423-439.\smallskip

\noindent Evans, M. and Moshonov, H. (2006) Checking for prior-data conflict.
Bayesian Analysis, 1, 4, 893-914.\smallskip

\noindent Evans, M. and Tomal, J. (2018) Multiple testing via relative belief
ratios. FACETS, 3: 563-583, doi: 10.1139/facets-2017-0121.\smallskip

\noindent Farr, C., Ruggeri F. and Mengersen, K. (2018) Prior and posterior
linear pooling for combining expert opinions: uses and impact on Bayesian
networks---the case of the wayfinding model. Entropy , 20, 209.
doi:10.3390/e20030209\smallskip

\noindent French, S. (2011) Aggregating expert judgement. Revista de la Real
Academia de Ciencias Exactas, F\'{\i}sicas y Naturales. Serie A,
Matem\'{a}ticas, 03, 105 (1), 181-206, Springer.\smallskip

\noindent Gelman, A. and O'Rourke, K. (2017) Attitudes towards amalgamating
evidence in statistics. Manuscript.\smallskip

\noindent Genest, C. (1984a) A characterization theorem for externally
Bayesian groups. Annals of Statistics, 12, 2, 1100-1105.\smallskip

\noindent Genest, C. (1984b) A conflict between two axioms for combining
subjective distributions. J. of the Royal Statistical Society. Series B, 46
(3), 403-405.\smallskip

\noindent Genest, C. (1984c) Pooling operators with the marginalization
property. Canadian Journal of Statistics, 12, 2, 153-163.\smallskip

\noindent Genest, C. and McConway, K. J. (1990) Allocating the weights in the
linear opinion pool. J. of Forecasting, 9, 53-73.\smallskip

\noindent Genest, C. and Zidek, J.V. (1986) Combining probability
distributions: a critique and an annotated bibliography. Statistical Science,
1, 114-135.\smallskip

\noindent Goodman, B. (2021) Daniel Kahneman says noise is wrecking your
judgment. Here's why, and what to do about it. Barrons, Economics \ Q\&A, May
28, 2021.\smallskip

\noindent Ladagga, R. (1977) Lehrer and the consensus proposal. Synthese, 36,
473-477.\smallskip

\noindent Lehrer, K. and Wagner, C. (1981) Rational Consensus in Science and
Society. D. Reidel Publishing Co.\smallskip

\noindent McConway, K. J. (1981) Marginalization and linear opinion pools. J.
of the American Stat. Assoc., 76:374, 410-414. doi:
10.1080/01621459.1981.10477661\smallskip

\noindent Nott,D., Wang, X., Evans, M., and Englert, B-G. (2020) Checking for
prior-data conflict using prior to posterior divergences. Statistical Science,
35, 2, 234-253.\smallskip

\noindent O'Hagan A., Buck C. E., Daneshkhah A., Eiser J. R., Garthwaite, P.
H., Jenkinson, D. J., Oakley, J.\ E. and Rakow, T. (2006) Uncertain
Judgements:Eliciting Experts' Probabilities. John Wiley \& Sons Ltd.\smallskip\ 

\noindent Royall, R. (1997) Statistical Evidence: A Likelihood Paradigm.
Monographs on Statistics and Applied Probability 71, CRC Press, Taylor \&
Francis Group.\smallskip

\noindent Stone, M. (1961) The opinion pool. Annals of Mathematical
Statistics, 32,1339-1342.\smallskip

\noindent Vieland, V.J. and Chang, H. (2019) No evidence amalgamation without
evidence measurement. Synthese, 196:3139--3161.
https://doi.org/10.1007/s11229-017-1666-7\smallskip

\noindent Winkler, R. L. (1968) The consensus of subjective probability
distributions. Management Science, 15, 2, B-61-B-75.\smallskip

\noindent Yager, R., Alajlan, N. (2015) An intelligent interactive approach to
group aggregation of subjective probabilities. Knowledge-based systems, 83
(1), 170-175.\smallskip

\noindent Zellner, A. (1996) An Introduction to Bayesian Inference in
Econometrics. Wiley Classics.

\end{document}